\documentclass[12pt]{article}
\usepackage{amsmath}
\usepackage{amsthm}
\usepackage{amsfonts}
\usepackage{float}
\usepackage[dvips,final]{graphicx}
\usepackage{dcolumn}
\usepackage{comment} 
\usepackage{xspace}
\usepackage{url}
\usepackage{upgreek}
\usepackage{multicol}
\usepackage{pict2e,multirow}
\usepackage{verbatim}
\usepackage{wrapfig}
\usepackage{url}
\usepackage{tikz}
\usepackage[colorlinks=true,citecolor=blue,urlcolor=blue]{hyperref}
\usepackage{subfigure}
\usepackage[T1]{fontenc}
\usepackage[hmargin=.720in]{geometry}

\newcommand{\be}{\begin{equation}}
\newcommand{\ee}{\end{equation}}

\newcommand{\Dx}{\Delta x}

\newcommand{\bu}{\mathbf{u}}

\newcommand{\Dt}{\Delta t}

\newcommand{\dt}{\Delta t}
\newcommand{\aij}{\alpha_{i,j}}
\newcommand{\bij}{\beta_{i,j}}

\newcommand{\m}[1]{\mathbf{#1}}
\newcommand{\mA}{\m{A}}

\newcommand{\mI}{\m{I}}

\renewcommand{\v}[1]{\boldsymbol{#1}}

\newcommand{\vb}{\v{b}}
\newcommand{\vc}{\v{c}}
\newcommand{\ve}{\v{e}}

\newcommand{\sspcoef}{\mathcal{C}}

\newcommand{\ceff}{\sspcoef_{\textup{eff}}}
\newcommand{\DtFE}{\Dt_{\textup{FE}}}

\newcommand{\mC}{\m{C}}

\renewcommand{\v}[1]{\mathbf{#1}}

\title{Optimal Explicit Strong Stability Preserving Runge--Kutta Methods with High Linear Order and optimal Nonlinear Order}

\author{%
Sigal Gottlieb\thanks{Mathematics Department, University of Massachusetts Dartmouth, 285 Old Westport Road,
North Dartmouth MA 02747. S.Gottlieb can be reached at sgottlieb@umassd.edu, Zachary Grant 
 at  zgrant@umassd.edu, and  Daniel Higgs at danielhiggs@gmail.com.},
Zachary Grant\footnotemark[1], and Daniel Higgs\footnotemark[1]
}
\begin{document}
\maketitle


\bibliographystyle{siam}

\begin{abstract} 

High order spatial discretizations with monotonicity properties are often desirable 
for the solution of hyperbolic PDEs. These methods can advantageously be coupled with
high order strong stability preserving time discretizations.
The search for high order strong stability time-stepping methods with large allowable strong stability coefficient
has been an active area of research over the last two decades.  This research has shown that
explicit SSP Runge--Kutta methods exist only up to fourth order. However, if we restrict ourselves
to solving only linear autonomous problems, the order conditions simplify and this order barrier
is lifted:  explicit SSP Runge--Kutta methods of any {\em linear order} exist.
These methods reduce to second order when applied to nonlinear problems. In the current work we 
 aim to  find explicit SSP Runge--Kutta methods with large allowable time-step, that feature high linear order 
and simultaneously have the optimal fourth order nonlinear order.
These methods have  strong stability coefficients that approach those of the linear methods as the 
number of stages and the linear order is increased. This work shows that when a high linear order method
is desired, it may be still be worthwhile to use methods with higher nonlinear order.

\end{abstract}

\section{Introduction\label{sec:intro}}

Explicit strong stability preserving (SSP) Runge--Kutta methods were developed for the time evolution of
hyperbolic conservation laws 
$ U_t +f(U)_x = 0$,
with discontinuous solutions \cite{shu1988, shu1988b}. 
These works studied total variation diminishing (TVD) spatial discretizations
that can handle discontinuities.  The spatial discretizations used to approximate $f(U)_x$ were 
carefully designed so that when the resulting system of ODEs
\begin{eqnarray}
\label{ode}
u_t = F(u),
\end{eqnarray}
(where $u$ is a vector of approximations to $U$,  $u_j \approx U(x_j) $) is  evolved in time using the 
{\em forward Euler  method}, the solution at time $u^n$ satisfies a strong stability property of the form
\begin{eqnarray} \label{monotonicity}
\| u^n + \dt F(u^{n}) \| \leq \| u^n \| 
\end{eqnarray}
under a step size restriction
\begin{eqnarray} \label{FEcond}
0 \leq \dt \leq \DtFE.
\end{eqnarray}
The term $\| \cdot \|$  can represent, as it did in  \cite{shu1988, shu1988b} 
the total variation semi-norm, or indeed any other semi-norm,
norm, or convex functional, as determined by the design of the spatial discretization. 
These spatial discretizations satisfy the strong stability property $\|u^{n+1} \| \leq \|u^{n} \| $
{\em when coupled with the forward Euler time discretization},
but in practice  a higher order time integrator, that will still satisfy this property, is desired. To accomplish this,
we attempt to re-write a higher order time discretization as a convex combination of forward Euler steps, 
so that any convex functional property that is satisfied by the forward Euler method will still be satisfied by the 
higher order time discretization.

An $s$-stage explicit Runge--Kutta method can be written in the form  \cite{shu1988},
\begin{eqnarray}
\label{rkSO}
u^{(0)} & =  & u^n, \nonumber \\
u^{(i)} & = & \sum_{j=0}^{i-1} \left( \aij u^{(j)} +
\dt \bij F(u^{(j)}) \right), \; \; \; \; i=1, . . ., s\\
 u^{n+1} & = & u^{(s)} . \nonumber
\end{eqnarray}
If all the coefficients $\aij$ and $\bij$ are non-negative, and a given $\aij$ is zero only if its corresponding $\bij$ is zero,
then each stage can be rearranged into a convex combination of forward Euler steps
\begin{eqnarray*}
\| u^{(i)}\| & =  & 
\| \sum_{j=0}^{i-1} \left( \aij u^{(j)} + \dt \bij F(u^{(j)}) \right) \|   \nonumber \\
& \leq &  \sum_{j=0}^{i-1} \aij  \, \left\| u^{(j)} + \dt \frac{\bij}{\aij} F(u^{(j}) \right\|  \\
& \leq & \|u^n\|  \; \; \; \; \forall \dt \leq \min_{i,j} \frac{\aij}{\bij} \DtFE,
\end{eqnarray*}
where  if any of the $\beta$'s are equal to zero, the corresponding 
ratio is considered infinite.
The last inequality above  follows from  the strong stability conditions \eqref{monotonicity} and  \eqref{FEcond} 
\[ \left\| u^{(j)} + \dt F(u^{(j})  \right\|   \leq \left\| u^{(j)}   \right\|   \; \; \; \forall \dt \leq \DtFE \]  and 
the consistency condition $\sum_{j=0}^{i-1} \aij =1$. From this we can conclude that 
whenever the explicit Runge--Kutta method can be decomposed into convex combinations
of forward Euler steps, then any convex functional  property  \eqref{FEcond} satisfied by forward Euler will be {\em preserved} by the higher-order time discretizations, perhaps under a different time-step restriction 
$\Dt \le \sspcoef \DtFE$ \cite{shu1988b}.  
Thus, this type of decomposition  where $\min_{i,j} \frac{\aij}{\bij}  >0$ 
is clearly a sufficient condition for strong stability preservation.
It has also been shown \cite{SSPbook2011, kraaijevanger1991,spijker2007} that this convex combination 
condition is necessary for strong stability preservation. 
If a method does not have a convex combination decomposition into forward Euler steps we can find 
some ODE with some initial condition such that the forward Euler condition is satisfied but the method
does not satisfy the strong stability condition for any positive time-step \cite{SSPbook2011}.
Methods that can be decomposed like this with with  $\sspcoef > 0$ are 
called strong stability preserving (SSP), 
and the coefficient $\sspcoef$ is known as the {\em SSP coefficient} of the method.

SSP methods guarantee the strong stability of the numerical solution for any ODE and any convex functional
provided only that the forward Euler condition \eqref{monotonicity}  is satisfied under a time step \eqref{FEcond}.
This is a very strong requirement that leads to severe restrictions on the allowable order of SSP methods,
and on the size of the allowable time step $\Dt \le \sspcoef \DtFE.$
We seek high order SSP Runge--Kutta methods  with the largest allowable time-step.
 The forward-Euler time step  $\DtFE$ is a property of the  spatial discretization method only, and so our aim in
 searching for time-stepping methods that preserve the strong stability property is to maximize the  
 {\em SSP coefficient} $\sspcoef$ of the method. 
A more relevant quantity may be the total cost of the time evolution,
 which in our case translates into the allowable time step  relative to the number of function evaluations
at each time-step (typically the number of stages of a method).
For this purpose we define the {\em effective SSP coefficient} $\ceff = \frac{\sspcoef}{s}$
where $s$ is the number of stages. This value allows us to compare the efficiency of 
explicit methods of a given order. It has been shown \cite{SSPbook2011} that all explicit $s$-stage Runge--Kutta 
methods have an SSP bound $\sspcoef \leq s$, and therefore $\ceff = 1 $,  but this upper bound is not always attained.

In   \cite{kraaijevanger1991,ruuth2001} it was shown that explicit SSP Runge--Kutta methods  
cannot exist for order $p>4$  However, in the special case where we consider only linear autonomous problems, 
explicit SSP Runge--Kutta methods exist for  any {\em linear} order $p_{lin}$ \cite{gottlieb2003}.
The linear and nonlinear order conditions are equivalent up to and including order $p=2$, so in 
 this work we consider explicit SSP Runge--Kutta methods that have nonlinear order $p=3$ and $p=4$,
and have higher linear orders $p_{lin}> p$. 
In Section \ref{sec:background} we review the SSP properties of  explicit Runge--Kutta 
methods and discuss the linear and nonlinear
order conditions. Using these order conditions and  optimization problem described in 
\cite{ketcheson2008, ketcheson2009a, ketcheson2009, tsrk}, in Section \ref{sec:optimization} we describe the 
optimization code in MATLAB (based on \cite{ketchcodes}) used to find explicit Runge--Kutta 
methods that have $p_{lin} > p$ with optimal SSP coefficient.  In Section \ref{sec:optimal} 
we list some of the new methods and their effective SSP coefficients, and
 in Section \ref{sec:test} we demonstrate  the performance of these methods on a selection of 
test problems.

\section{A review of explicit SSP Runge--Kutta methods} \label{sec:background}
Strong stability preserving methods were first developed by Shu \cite{shu1988, shu1988b} for use with total variation diminishing
spatial discretizations. In these works, the authors presented second and third order methods that have 
$\sspcoef=1$.
The explicit SSP Runge--Kutta method of order $p=s=2$ 
\begin{eqnarray*}
u^{(0)} & =& u^{n} \\
u^{(1)} & = &u^{(0)} + \Dt F(u^{(0)}) \\
u^{n+1} & = &\frac{1}{2} u^{(0)} + \frac{1}{2} \left( u^{(1)} +  \Dt F(u^{(1)}) \right). 
\end{eqnarray*}
and the $p=s=3$ method
\begin{eqnarray*}
     u^{(1)} &= & u^n + \dt F(u^n) \nonumber \\
     u^{(2)} &= & \frac{3}{4} u^n + \frac{1}{4} u^{(1)} + \frac{1}{4} \dt F(u^{(1)}) \\
     u^{n+1} & = & \frac{1}{3} u^n + \frac{2}{3} u^{(2)} +  \frac{2}{3} \dt F(u^{(2)}).
\end{eqnarray*}
These methods  were proven optimal in \cite{gottliebshu1998}. 

It was shown in \cite{kraaijevanger1991,ruuth2001} that no four stage fourth order explicit Runge--Kutta
methods exist with positive SSP coefficient. By considering methods with $s>p$, fourth order methods
with order $p=4$ have been found. Notable among these is the $(s,p)=(5,4)$ method with $\sspcoef=1.508$
($\ceff=0.302$) in \cite{SpiteriRuuth2002}
\begin{eqnarray*}
u^{(1)} & = & u^n +  0.391752226571890 \dt F(u^n) \\ 
u^{(2)} & = &  0.444370493651235 u^n +  0.555629506348765 u^{(1)} 
+ 0.368410593050371 \dt F(u^{(1)}) \\ 
u^{(3)} & = &  0.620101851488403 u^n +  0.379898148511597 u^{(2)} 
 + 0.251891774271694  \dt F(u^{(2)}) \\ 
u^{(4)} & = &  0.178079954393132 u^n + 0.821920045606868 u^{(3)} 
+  0.544974750228521 \dt F(u^{(3)})\\ 
u^{n+1} & = &    0.517231671970585 u^{(2)} 
 +  0.096059710526147 u^{(3)} +  0.063692468666290 \dt F(u^{(3)}) \\ 
& & +  0.386708617503269 u^{(4)} +   0.226007483236906 \dt F(u^{(4)}) \, ,
\end{eqnarray*}
and the $(s,p)=(10,4)$ method with $\sspcoef=6$
($\ceff=0.6$) in \cite{ketcheson2008}
\begin{eqnarray*}
u^{(1)} & = & u^n + \frac{1}{6} \dt F(u^n) \\ 
u^{(i+1)} & = & u^{(i)} + \frac{1}{6} \dt F(u^{(i)}) \; \; \;  i=1,2, 3\\ 
u^{(5)} & = & \frac{3}{5} u^n +   \frac{2}{5} u^{(4)} +  \frac{1}{15} \dt F(u^{(4)}) \\ 
u^{(i+1)} & = & u^{(i)} + \frac{1}{6} \dt F(u^{(i)}) \; \; \; i=5,6,7,8 \\ 
u^{n+1} & = &   \frac{1}{25} u^{n} +  \frac{9}{25} u^{(4)} +  \frac{3}{5} u^{(9)}
+ \frac{3}{50}  \dt F(u^{(4)}) + \frac{1}{10}  \dt F(u^{(9)})  \, .
\end{eqnarray*}
It was shown  \cite{kraaijevanger1991,ruuth2001} 
that no methods of order $p \geq 5$ with positive SSP coefficients can exist. 

This restrictive order barriers on explicit SSP Runge--Kutta methods stem in part from 
the nonlinearity of the ODEs.  For order of accuracy on linear autonomous ODE systems, 
explicit SSP Runge--Kutta methods  need only satisfy a smaller set of  order conditions.
If we require only that the method   have high linear order ($p_{lin}$), then the order barrier is 
broken and explicit Runge--Kutta methods with  positive SSP coefficients exist for arbitrarily high linear orders.
Optimally contractive explicit Runge--Kutta methods were studied 
by Kraaijevanger in \cite{kraaijevanger1986}, where he gives optimal {\em linear} methods 
for many values of $s$ and $p$, including $1\le p_{lin} \le s \le 10$, and
$p_{lin}  \in \{1,2,3,4,s-1,s-2,s-3,s-4\}$ for any $s$.  
These methods are interesting because their SSP coefficients
serve as upper bounds for nonlinear methods, but they may also be useful in their own right.
Although SSP methods were first developed for nonlinear problems,  the 
strong stability preserving property can be useful for linear problems such as
Maxwell's equations and linear elasticity.

First and second methods that have $s$ stages have been shown to attain the theoretical bound
$\sspcoef \leq s-p_{lin}+1$. Methods of order $p_{lin}$ with $s=p_{lin}$ and $s=p_{lin} +1$
also exist with $\sspcoef \leq s-p_{lin}+1$.
These methods can be  found in \cite{kraaijevanger1986,gottlieb2003, SSPbook2011}, 
and are given here  in their  canonical Shu-Osher form:
The family of $s$-stage, linear order $p_{lin}=s-1$ methods 
has $\sspcoef=2$ and $\ceff=\frac{2}{s}$:
\begin{eqnarray*}
u^{(0)} &=& u^n  \\ 
u^{(i)} & = & u^{(i-1)} + \frac{1}{2} \dt F \left(  u^{(i-1)}\right), \qquad i=1, . . . , m-1 \\ 
u^{(s)} & = &  \sum_{k=0}^{s-2}  \alpha^s_{k} u^{(k)} +
\alpha^s_{s-1}  \left( u^{(s-1)}+ \frac{1}{2} \dt F \left( u^{(s-1)} \right) \right), \\
 \bu^{n+1} & = & u^{(s)}  ,
\end{eqnarray*}
where the coefficients $\alpha^s_k$  of the final stage of the $s$-stage method are
given iteratively by
\[ \alpha^s_{k} =  \frac{2}{k} \alpha^{s-1}_{k-1} \; \; \; \mbox{for} \; \; \; k=1, . . . , s-2, \; \; \; \;
\alpha^{s}_{s-1}  =  \frac{2}{s} \alpha^{s-1}_{s-2} , \qquad
\alpha^s_{0}  =  1- \sum_{k=1}^{s-1} \alpha^s_{k}, \]
starting from the coefficients of the $2$-stage, first order method
$\alpha^2_{0} = 0 $ and $\alpha^2_{1}=1$.

The family of $s$-stage, linear order $p_{lin}=s$ methods 
has $\sspcoef=1$ and $\ceff=\frac{1}{s}$:
\begin{eqnarray*}
u^{(0)} &=& u^n  \\ 
u^{(i)} & = & u^{(i-1)} + \dt F \left( u^{(i-1)} \right), \qquad
i=1, . . . , s-1 \\ 
u^{(s)} & = &  \sum_{k=0}^{s-2}  \alpha^s_{k} u^{(k)} +
\alpha^s_{s-1}  \left(u^{(s-1)}+\dt F \left( u^{(s-1)} \right)\right), \\
 u^{n+1} & = & u^{(s)}  .
\end{eqnarray*}
Here the coefficients $\alpha^s_k$  of the final stage of the $s$-stage method are
given iteratively by
\[ \alpha^s_{k} = \frac{1}{k} \alpha^{s-1}_{k-1}  \; \;  \; \mbox{for} \; \;  k=1, . . . , s-2 ,   \; \; \; \; \; 
\nonumber \alpha^m_{s-1}  =  \frac{1}{s!}, \qquad
\alpha^s_{0}  =  1- \sum_{k=1}^{s-1} \alpha^s_{k},\]
starting from the coefficient of the forward Euler method
$\alpha^1_{0} = 1$. 

However, all these methods with high linear order $p_{lin}$ have low nonlinear order $p=2$.
The idea that we pursue in this paper is the construction of explicit SSP Runge--Kutta methods that
have a high linear order $p_{lin} > 4$ while retaining the highest possible nonlinear order
$p=4$. We also consider methods with $p=3$ and $p_{lin}  \geq  4$ for comparison.
The idea behind these methods is that they would be the best possible methods (in terms of SSP coefficient
and order) for linear problems, without compromising order when applied to nonlinear problems.

\section{Formulating the optimization problem} \label{sec:optimization}

The Shu-Osher  form of an explicit Runge--Kutta method, given in \eqref{rkSO}, is most convenient for 
observing the SSP coefficient. However, this form is not unique, and not the most efficient form to use for the 
optimization procedure \cite{ketcheson2008}.  The Butcher form of the explicit method given by
\begin{eqnarray} \label{butcher}
u^{(i)} & = & u^{n} + \dt \sum_{j=1}^{i-1} a_{ij} F(u^{(j)}) \; \; \; \; (1 \leq i \leq s) 
\\ \nonumber
u^{n+1} & = & u^n + \dt \sum_{j=1}^s b_j  F(u^{(j)}).
\end{eqnarray}
(where the coefficients $a_{ij}$ are place into the matrix $\mA$ and $b_{j}$ into the row vector $\vb$)
is unique, so rather than perform a search for the optimal convex combination of the Shu-Osher form \eqref{rkSO},
we define the optimization problem in terms of the Butcher coefficients.
The conversion from the Shu-Osher form to Butcher form, and from an optimal Butcher form to the 
canonical Shu-Osher form is discussed in  \cite{SSPbook2011}.

 We follow the approach developed by David Ketcheson and  successfully used in 
\cite{ketcheson2008, ketcheson2009a, ketcheson2009, Ketcheson2010, SSPbook2011, tsrk}: 
we search for coefficients $\mA$ and $\vb$ that maximize
the value $r$ subject to constraints:
\begin{eqnarray}
(1) & \left( \begin{array}{ll} \mA & 0 \\ \vb & 0  \\ \end{array} \right)
\left(\mI + r 
\left( \begin{array}{ll} \mA & 0 \\ \vb & 0  \\ \end{array} \right)
\right)^{-1} \geq 0   \\ 
&\mbox{where the inequality is understood component wise.} \nonumber \\
(2) & \left\| r  \left( \begin{array}{ll} \mA & 0 \\ \vb & 0  \\ \end{array} \right)
\left(\mI + r  \left( \begin{array}{ll} \mA & 0 \\ \vb & 0  \\ \end{array} \right)
\right)^{-1} \right\|_\infty \leq 1 \\
 (3) & \tau_k(\mA, b) = 0 \; \; \; \mbox{for} \; \; \;  k=1, . . ., P, 
 \end{eqnarray}
 where $\tau_k$ are the order conditions.
After this optimization we have the coefficients $\mA$ and $\vb$ and an optimal value $\sspcoef = r$
that define the method. 

\subsection{Linear and nonlinear order conditions} \label{orderconditions}
The order conditions  $\tau_k(\mA, \vb)$ appear as  the equality  constraints on the
optimization problem.  In this work, we consider methods 
that have $p=3$ and $p=4$ but have higher linear order $p_{lin} > 4$. 
 In this subsection, we list these order conditions. 

\noindent{\bf Linear Order Conditions:}
Given a Runge--Kutta method written in the Butcher form with coefficients $\mA$ and $\vb$ (and $\vc= \mA \ve$ 
where $\ve$ is a vector of ones), 
the order conditions that guarantee $p_{lin}$ order accuracy for a linear problem can be simply expressed as
\begin{eqnarray}
\vb^T \mA^{q-2} \vc = \vb^T \mA^{q-1} \ve = \frac{1}{q!} \; \; \; \; \; \forall q=1, . . . , p_{lin}.
\end{eqnarray}

\noindent{\bf Nonlinear Order Conditions:}
If we want a method to demonstrate the correct order of accuracy for nonlinear problems,
the  first  and second order conditions are the same as above:
\begin{align*}
\vb^T \ve & = 1 & \vb^T \vc & = \vb^T \mA \ve =  \frac{1}{2}.
\end{align*}
A method that satisfies these conditions will be second order for both linear and nonlinear problems.
Two  additional conditions are required for third order accuracy\footnote{These
nonlinear order conditions  follow Albrecht's notation as Ketcheson found these to be handled more efficiently by the optimizer.}:
\begin{align} \label{thirdorder}
\vb^T \vc^2 & = \frac{1 }{3}, &
\vb^T  \left( \frac{\vc^2}{2!} - \mA \vc \right) = 0.
\end{align}
Note that when the first of these conditions is satisfied, the second condition is equivalent to 
$\vb^T \mA \vc = \frac{1}{3!}$, which is the linear third order condition.
Four more conditions are required for the method to be fourth order for a nonlinear problem  
\begin{align} \label{fourthorder}
\vb^T \vc^3  = \frac{1}{4}, \; \; \; 
\vb^T \mA \left( \frac{\vc^2}{2!} - \mA \vc \right)  = 0, \; \; \; 
\vb^T \left( \frac{\vc^3}{3!} - \frac{\mA \vc^2}{2!}  \right)  = 0, \\
\vb^T \mC \left( \frac{\vc^2}{2!} - \mA \vc \right)  = 0 \; \; \; (\mbox{where} \; \; \mC=diag(\vc)).  \nonumber
\end{align}
Note that the first three conditions  together imply the fourth order linear
order condition $b^T A^2 c = \frac{1}{4!}.$

In this work we consider the nonlinear order conditions only up to $p=4$ because  it is known that 
there are no explicit SSP Runge--Kutta methods greater than fourth order, 
but we consider higher order linear order conditions $p_{lin} > 4$.

\vspace{-.25in}
\section{Optimal methods\label{sec:optimal}}

Using David Ketcheson's MATLAB optimization code \cite{ketchcodes} with our modified order conditions  (described in 
Section \ref{orderconditions})  we produce the optimal linear/nonlinear (LNL) methods in this section. 
This \begin{wraptable}[13]{r}{0.625\textwidth} \vspace{-.1in} 
{\footnotesize \begin{tabular}{|l|l|l|l|l|l|l|l|l|}
\hline
$s$   &$p_{lin}$ = 5&6&7&8&9&10&11&12\\\hline
2&--&--&--&--&--&--&--&--\\\hline
3&--&--&--&--&--&--&--&--\\\hline
4&--&--&--&--&--&--&--&--\\\hline
5&1&--&--&--&--&--&--&--\\\hline
6&2&1&--&--&--&--&--&-- \\\hline
7&2.6506&2&1&--&--&--&--&--\\\hline
8&3.3733&2.6506&2&1&--&--&--&--\\\hline
9&4.1&3.3733&2.6506&2&1&--&--&--\\\hline
10&4.8308&4.1&3.3733&2.6506&2&1&--&--\\\hline
11&5.5193&4.8308&4.1&3.3733&2.6506&2&1&--\\\hline
12&6.349&5.5193&4.686&4.1&3.3733&2.6506&2&1\\\hline
\end{tabular}}
\caption{\vspace{-.15in} SSP coefficients for linear and $p=3$ methods.}
\label{tab:p2p3}
\end{wraptable}
code formulates the optimization problem in Section \ref{sec:optimization} in MATLAB and uses {\tt fmincon}
to find the coefficients $\mA$ and $\vb$ that yield the largest possible $\sspcoef$. We set the tolerances
on {\tt fmincon} to $10^{-14}$.
We used this code to generate methods with $p=3,4$ and $p_{lin} = 5, . . ., 12$. We compare these methods
with $p=2$ "linear" methods that we generated and matched to known optimal methods. 
Our primary interest is the size of the SSP coefficient for each method. We denote the SSP coefficient
for a method with $s$ stages, linear order $p_{lin}$ and nonlinear order $p)$ method $\sspcoef_{(s,p_{lin}, p)}$.

The SSP coefficients for the $p=3$  methods with  a given number of stages and 
linear order   are the same as for the corresponding linear methods, (i.e. 
$\sspcoef_{(s,p_{lin}, 2)} = \sspcoef_{(s,p_{lin}, 3)}$). This indicates that, for these values of $s$ and
$p_{lin}$ the additional condition $\vb^T \vc^2 = \frac{1}{3}$ needed for nonlinear third order 
does not pose additional constraints  on the strong stability properties of the method. 
Table  \ref{tab:p2p3} shows the SSP coefficients of the $p=2$ and $p=3$ methods
with $5 \leq p_{lin} \leq s \leq 12$. The coefficients for $s \leq 10$ are 
known to be optimal because they match the linear threshold in Kraiijevanger's paper \cite{kraaijevanger1986}.

\begin{figure}[ht]
\begin{minipage}[b]{0.45\linewidth}
\centering 
\includegraphics[width=1.13\textwidth]{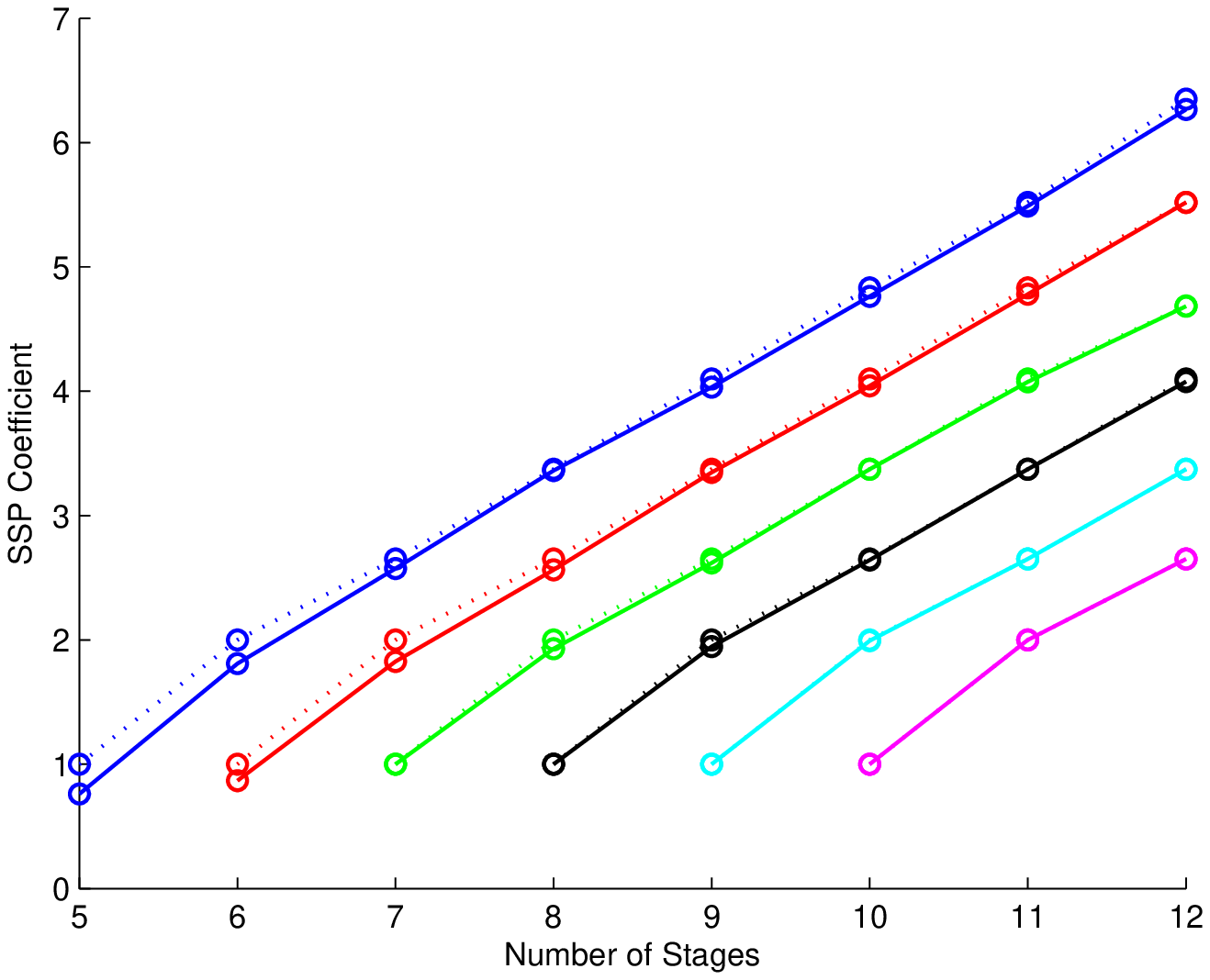}
\caption{The SSP coefficient as a function of the number of coefficients 
for $p=2$ (dotted line) and $p=4$ (solid line) for $p_{lin}=5$ (blue),
$p_{lin}=6$ (red), $p_{lin}=7$ (green), $p_{lin}=8$ (black), $p_{lin}=9$ (cyan),
and $p_{lin}=10$ (magenta).
As we increase the number of stages the SSP coefficients for the $p=4$ methods
approach those of the $p=2$ methods
\vspace{.01in}}
\label{fig:p2p4a}
\end{minipage}
\hspace{0.5cm}
\begin{minipage}[b]{0.45\linewidth}
\centering
\includegraphics[width=1.13\textwidth]{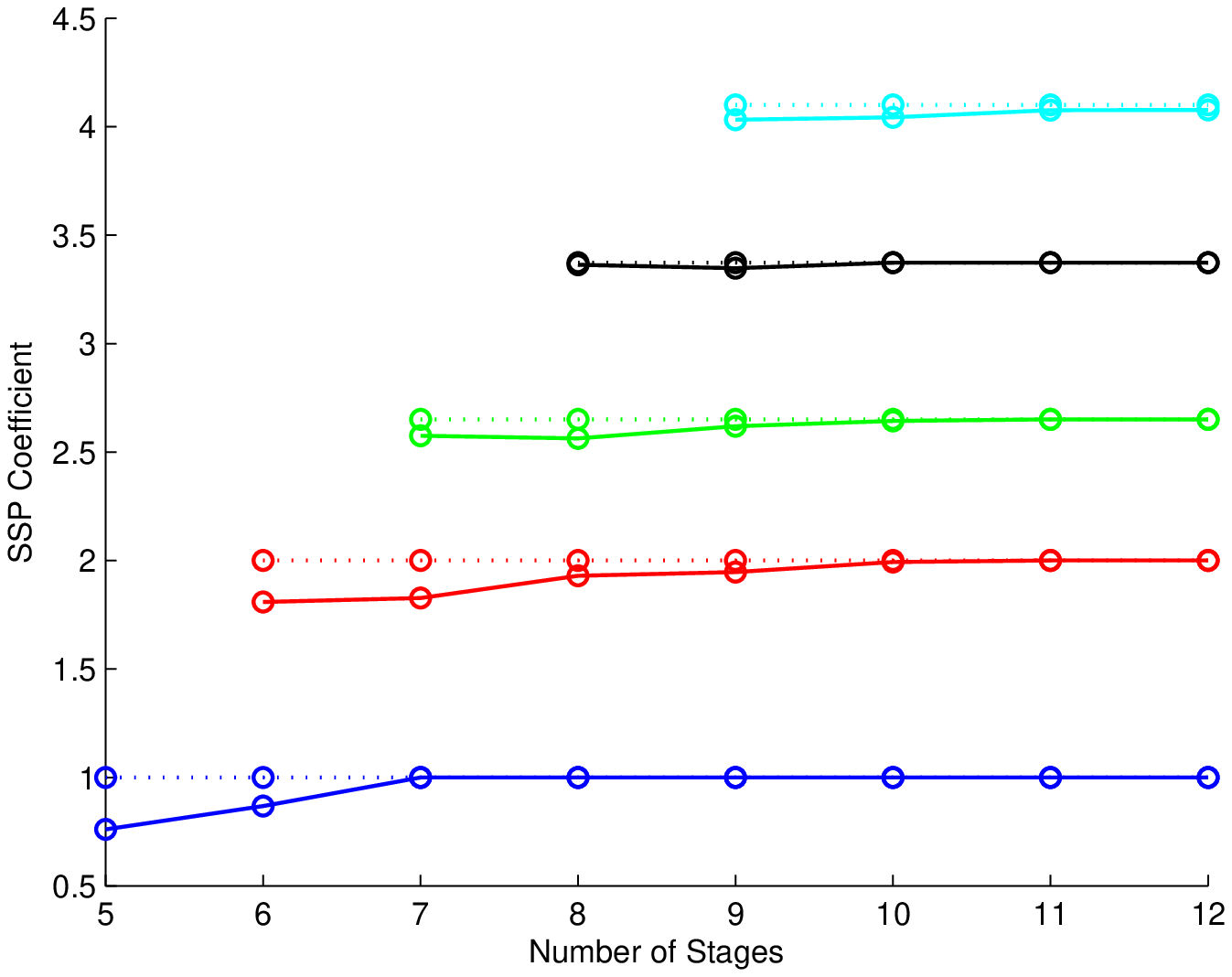}
\caption{The SSP coefficient of  linear methods $p=2$ (dotted line) and $p=4$ (solid line)
for $p_{lin} = s$ (blue), $p_{lin} = s-1$  (red), 
$p_{lin} = s-2$  (green), $p_{lin} = s-3$ (black), and $p_{lin} = s-4$ (cyan).
The SSP coefficient of the  $p=4$ methods (solid lines) approach those of the corresponding $p=2$ 
method as the number of stages increases
\vspace{.01in}}
\label{fig:p2p4b}
\end{minipage}
\end{figure}

\begin{wraptable}[13]{r}{0.67\textwidth} \vspace{-.12in} 
{\footnotesize 
\begin{tabular}{|l|l|l|l|l|l|l|l|l|}
\hline
1&5&6&7&8&9&10&11&12\\\hline
2&--&--&--&--&--&--&--&--\\\hline
3&--&--&--&--&--&--&--&--\\\hline
4&--&--&--&--&--&--&--&--\\\hline
5&0.76026&--&--&--&--&--&--&--\\\hline
6&1.8091&0.86773&--&--&--&--&--&--\\\hline
7&2.5753&1.8269&{\bf 1}&--&--&--&--&--\\\hline
8&3.3627&2.5629&1.9293&{\bf 1}&--&--&--&--\\\hline
9&4.0322&3.347&2.6192&1.9463&{\bf 1}&--&--&--\\\hline
10&4.7629&4.0431&{\bf 3.3733}&2.6432&1.9931&{\bf 1}&--&--\\\hline
11&5.4894&4.7803&4.0763&{\bf 3.3733}&{\bf 2.6506}&{\bf 2}&{\bf 1}&--\\\hline
12&6.267&{\bf 5.5193}&4.6842&4.0766&{\bf 3.3733}&{\bf 2.6506}&{\bf 2}&{\bf 1}\\\hline
\end{tabular} }
\caption{SSP coefficients for  $p=4$ methods.}
\label{tab:p4} 
\end{wraptable}
Table \ref{tab:p4} shows the SSP coefficients of the $p=4$ methods for $5 \leq p_{lin} \leq s \leq 12$.
In bold are the coefficients that match those of the $p=2$ methods.
In general, as we increase the number of stages the SSP coefficients for the $p=4$ methods
approach those of the $p=2$ methods, as shown in Figure \ref{fig:p2p4a}.
 
The tables clearly show that  the size of the SSP coefficient  depends on
the relationship between $p_{lin}$ and $s$, so it is illuminating to look at  the methods 
along the diagonals of these tables. Clearly, for $s=p_{lin}$ methods we have an optimal
value of $\sspcoef=1$ and $\ceff=\frac{1}{s}$. The $p=2$ and $p=3$ methods all attain this
optimal value, but for $p=4$ we have $\sspcoef_{(5,5,4)} = 0.76$ and $\sspcoef_{(6,6,4)} = 0.87$.
However, once we get to a high enough number of stages, all the methods  
with $p=4$ and $s=p_{lin} = 7, . . . , 10$ that have $\sspcoef=1$ and $\ceff=\frac{1}{s}$. 
Figure \ref{fig:p2p4b} shows that for the linear methods ($p=2$) the SSP coefficient is fixed for 
$p_{lin} = s$ (blue dotted line), $p_{lin} = s-1$  (red dotted line), 
$p_{lin} = s-2$  (green dotted line), $p_{lin} = s-3$ (black dotted line), and $p_{lin} = s-4$ (cyan dotted line), 
and that the SSP coefficient of the corresponding $p=4$ methods (solid lines) approach these as 
the number of stages increases.

%



\begin{figure}[H]
\includegraphics[scale=.45]{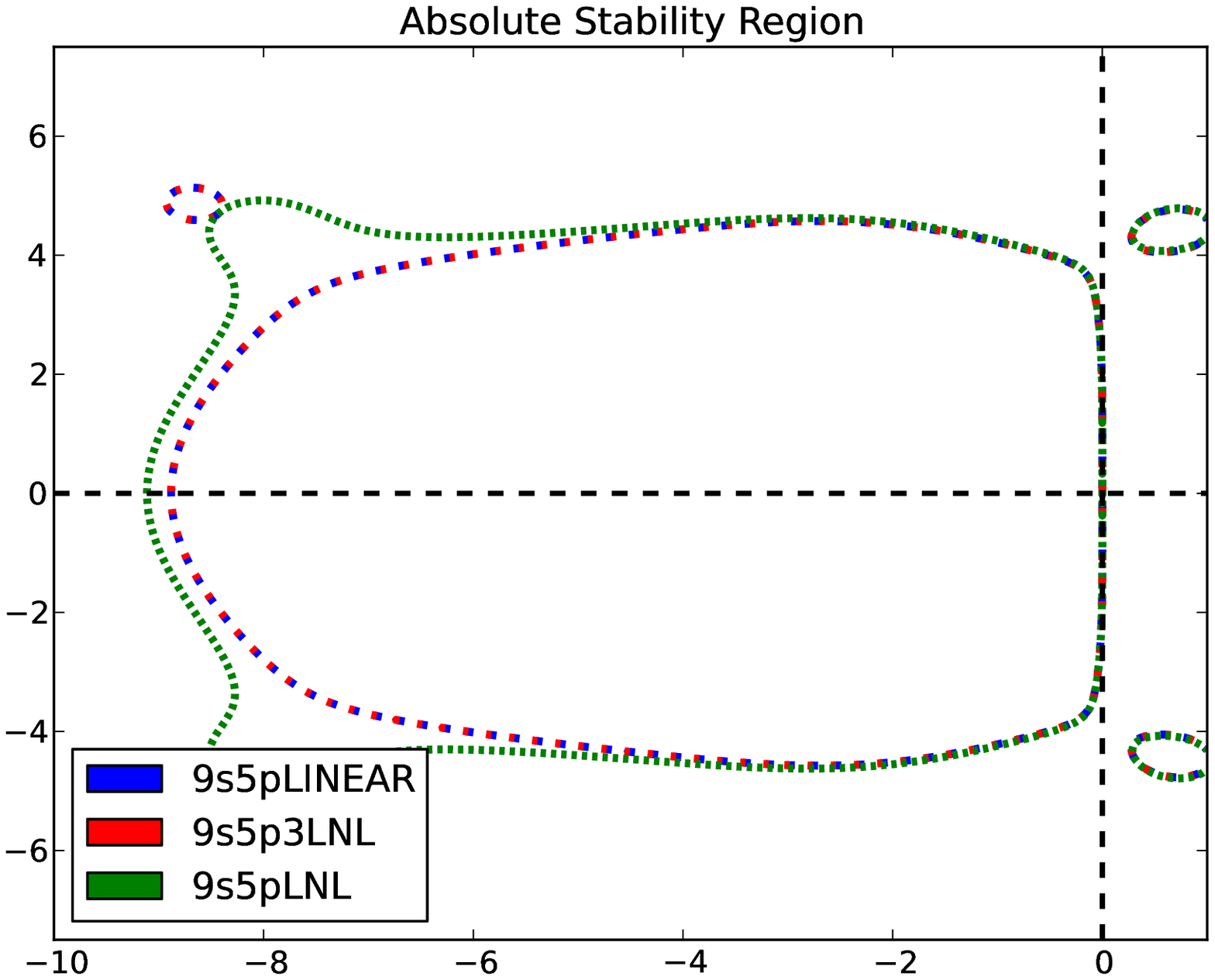} 
\includegraphics[scale=.45]{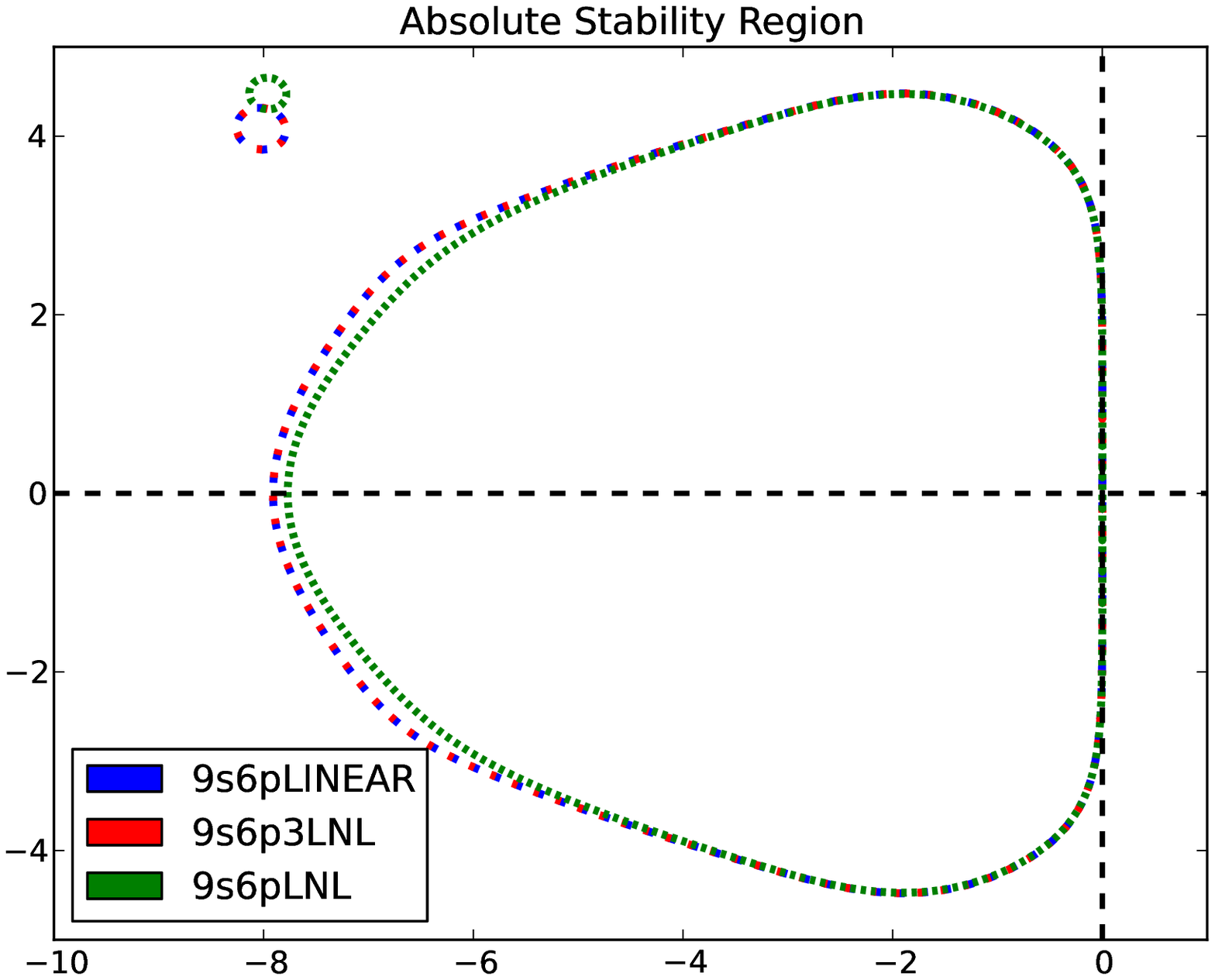}\\
\includegraphics[scale=.45]{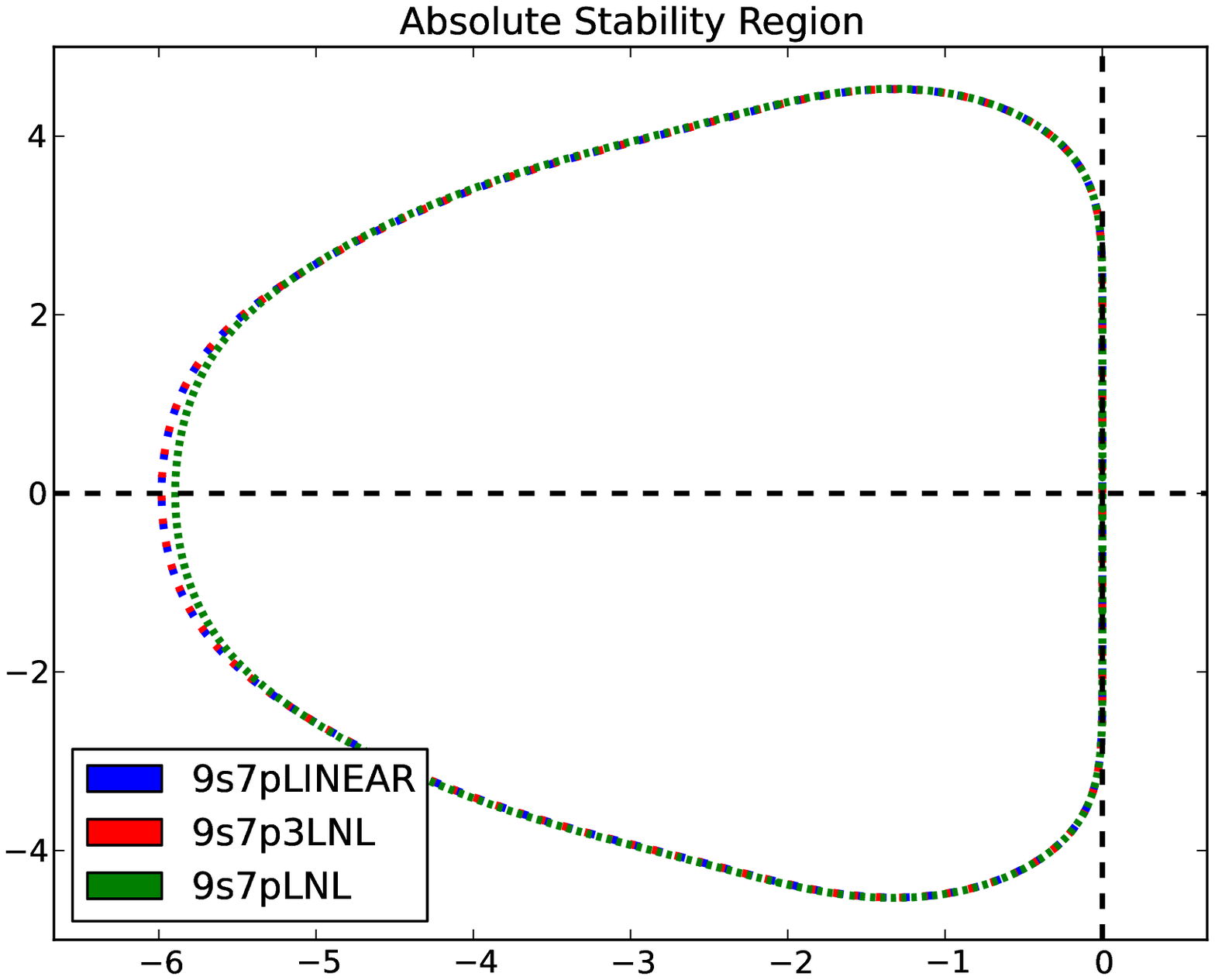} 
\includegraphics[scale=.45]{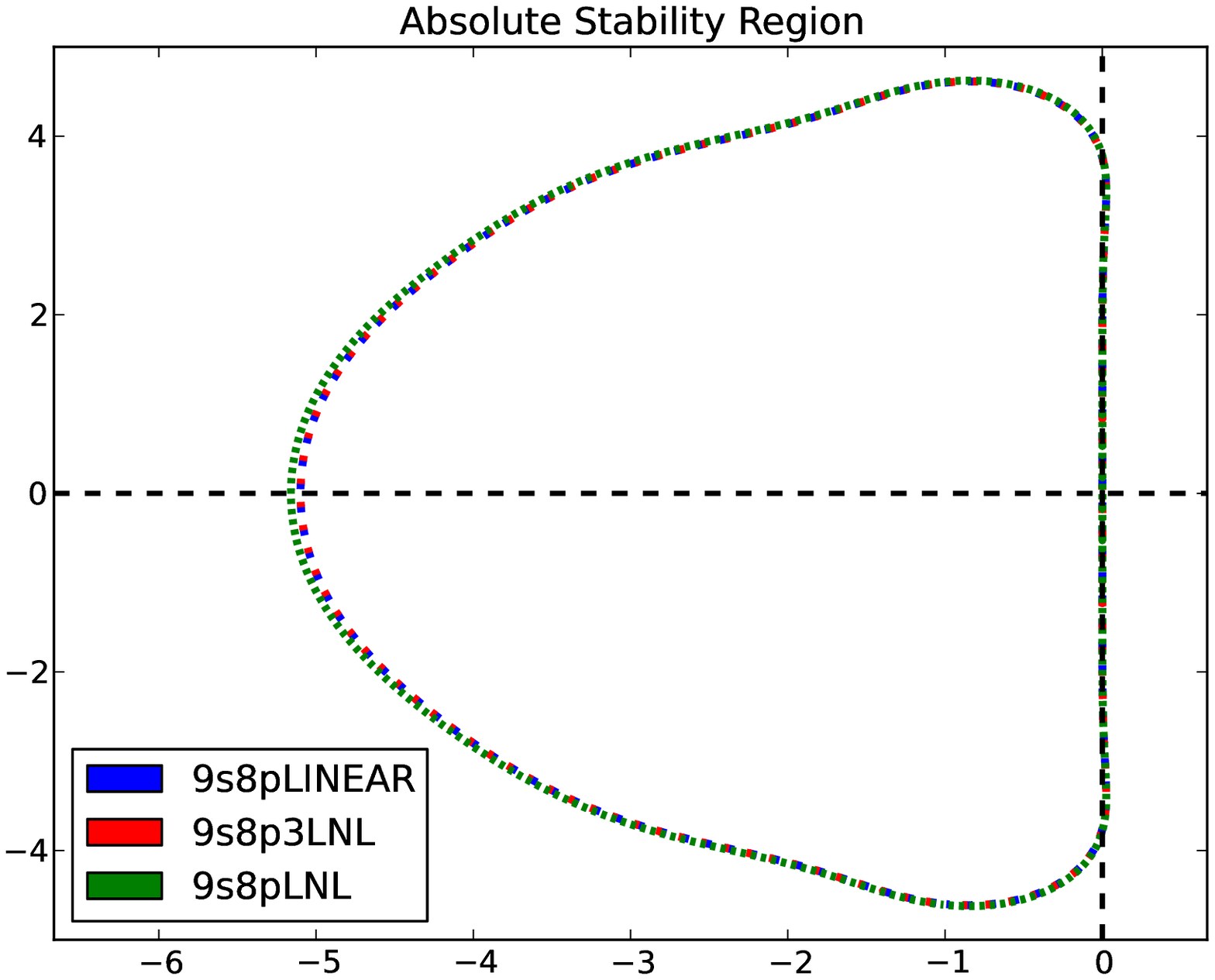} \\
\caption{Linear stability regions of nine stage methods of linear orders $p_{lin} = 5,6,7,8$ for linear (blue), 
$p=3$ (red) and $p=4$ (green) methods.
The $p=4$ methods approach the $p=3$ and $p=2$ methods as $p_{lin}$ increases.}
 \label{fig:LinStab2}
\end{figure}

It is interesting to note that the linear stability regions of the $p=2$, $p=3$ methods are generally identical.
The $p=4$ methods have stability regions that are virtually identical to those of the linear methods
when the SSP coefficient is identical. In addition, methods with $s=p_{lin}$ and $p=4$ all have the same stability 
regions as the corresponding linear methods, which is not surprising as the stability polynomial
of  $s=p_{lin}$ is unique. For the rest of the methods, we observe that for a given number of stages $s$,
as the  linear order $p_{lin}$ increases the  linear stability  regions of the $p=4$ methods
look closer to those of the linear methods. 
A nice illustration of this is the family of $s=9$   methods, shown in Figure \ref{fig:LinStab2}.


It is known in the literature that some methods with nonlinear orders $p=3$ and $p=4$ achieve the
linear threshold value. A nice example of this is Ketcheson's SSP Runge--Kutta method of $s=10$,
$p=p_{lin}=4$, which achieves the threshold value $\sspcoef = 6.0$. This suggests that the linear
order conditions are very significant to the value of the SSP coefficient. Indeed, we see this relationship in
Tables \ref{tab:p2p3} and \ref{tab:p4}, as we move right from column to column we see a significant drop
in SSP coefficient. 
For each application, one must decide if a higher linear order  is valuable,  
as we pay a price for requiring additional $p_{lin}$. However, once one has decided that the 
cost of a higher linear order is useful, there is no penalty in terms of SSP coefficient 
for requiring a higher nonlinear order $p=3$ and, in most cases,  little reason not to use $p=4$.
Our results show  that if one wishes to  use
 a method with high linear order $p_{lin} > 4$, then requiring $p=3$ or even $p=4$ rather than the standard $p=2$
is not usually associated with significant restriction on the SSP coefficient.
This can be beneficial in cases where the solution has linear and nonlinear components that need
to be accurately captured simultaneously, or in different regions, or at different time-levels, 
so that  the use of an SSP method
that has optimal nonlinear order and higher linear order would be best suited for all components
of the solution.

\section{Numerical Results\label{sec:test}}
In this section, the optimized LNL methods described in Section \ref{sec:optimal} are tested for convergence and SSP properties.
First, we test these methods for convergence on  both ODEs and PDEs to confirm that the 
new methods exhibit the desired linear and nonlinear orders. Next, we study the behavior of these methods
in conjunction with a higher order WENO spatial discretizations, and show that although the WENO method
is nonlinear, when applied to a linear smooth problem the higher order linear order $p_{lin}$ is beneficial. 
On the other hand, for nonlinear problems, both with shocks and without, the higher order nonlinear 
order $p=4$ is advantageous. Finally, the LNL methods are tested on  linear and nonlinear problems with spatial
discretizations that  are provably total variation diminishing (TVD) and positivity preserving. We show that for the linear
case, the observed time-step for the time stepping method to preserve the TVD property is well-predicted by the theoretical
SSP coefficient, while for positivity and for the nonlinear problem the theoretical time-step serves as a lower bound.

\begin{figure}[ht]
\begin{minipage}[b]{0.45\linewidth}
\centering 
\includegraphics[width=1.13\textwidth]{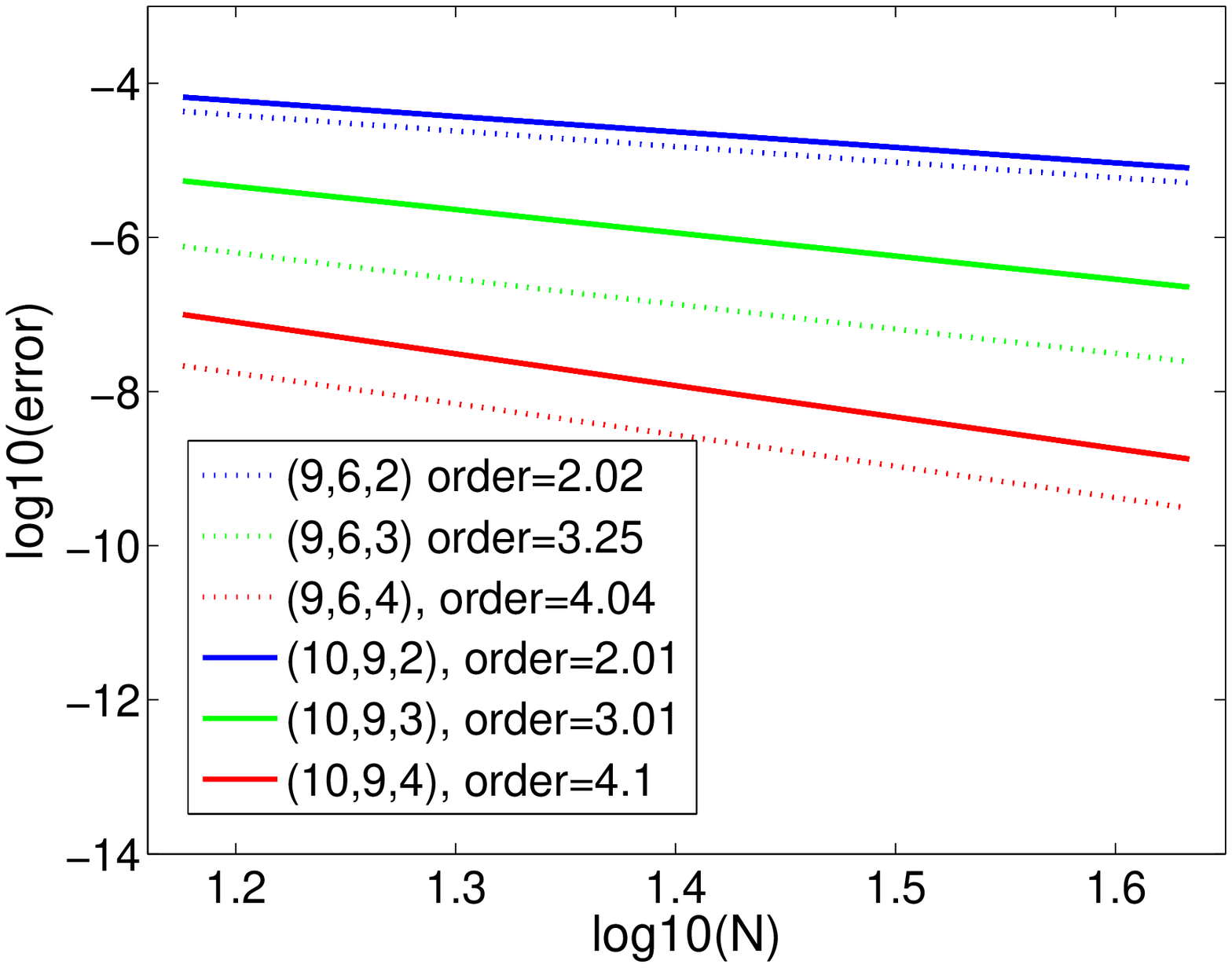}
\caption{Order verification of LNL Runge--Kutta methods on the van der Pol problem.
\vspace{.01in}}
\label{fig:VDPconv}
\end{minipage}
\hspace{0.5cm}
\begin{minipage}[b]{0.45\linewidth}
\centering
\includegraphics[width=\textwidth]{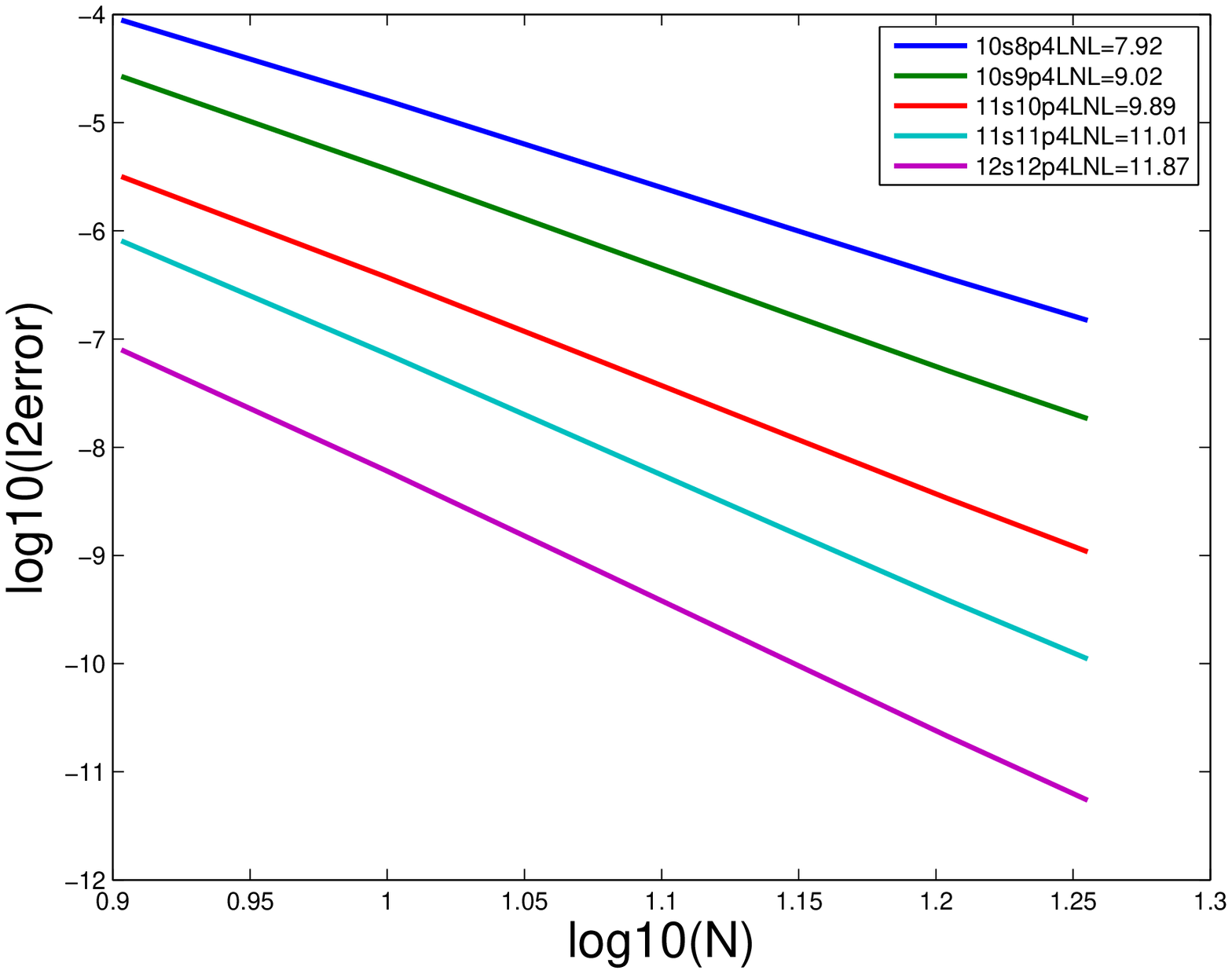}
\caption{Order verification of LNL Runge--Kutta methods on a linear advection problem with
Fourier spectral method in space.}
\label{fig:spectral}
\end{minipage}
\end{figure}


\noindent{\bf Example 1: Nonlinear ODE convergence study.}  
The van der Pol problem is a nonlinear system of ODEs:
\begin{eqnarray}
&u_1' = u_2 \\ 
&u_2' = \frac{1}{\epsilon} (-u_1 + (1-u_1^2) u_2)
\end{eqnarray}
We use $\epsilon = 10$ and initial conditions $u_0 = (0.5;  0)$. 
This was run  using $(s,p_{lin}, p) =$ (9,6,2), (9,6,3), (9,6,4), (10,9,2), (10,9,3), and (10,9,4) 
LNL Runge--Kutta methods to final time $T_{final} =  4.0$, with $\Delta t = \frac{T_{final}}{N-1} $ where 
$N=15,19,23,27,31,35,39,43$. 
The exact solution (for error calculation)  was calculated by  MATLAB's 
ODE45 routine with tolerances set to {\tt AbsTol=}$10^{-14}$
and {\tt  RelTol=}$10^{-13}$.
In Figure \ref{fig:VDPconv}  we show that the $log_{10}$ of the errors
in the first component vs. the $log_{10}$ of the number of points. The slopes of these lines (i.e the orders) are 
calculated by MATLAB's {\tt polyfit} function.
As expected, the rate of convergence follows the nonlinear order of the method.
In fact, we observe that a higher linear order is of no benefit at all for this example.

\noindent{\bf Example 2: PDE convergence study.} 
In this study we solve the  linear advection equation with sine wave  initial conditions and  periodic boundaries
\begin{eqnarray}
u_t & = & - u_x  \; \; \; x \in [0,1]  \\ \nonumber
u(0,x)& =& \sin(4\pi x) \; \; \; \; u(t,0)=u(t,1)  \nonumber
\end{eqnarray}
The Fourier spectral method was used  to discretize in space using $N=(9,11,13, 15, 17, 19)$ points.
The exact solution to this problem is a sine wave with period $4$ that travels in time, so  
the Fourier spectral method gives us an exact solution in space \cite{HGG2007} 
once we have two points per wavelength, allowing us to isolate the effect of the  
time discretization on the error. We run this problem for  five methods with orders
$p_{lin} = 8, 9, 10, 11, 12$,
$(s,p_{lin},p) = (10,8,4), (10,9,4), (11,10,4), (11,11,4), (12,12, 4)$.
Our final time is $T_{final}=1$, and the time step $\Delta t = 0.9 \Delta x$, where $\Delta x = \frac{1}{N-1}$. 
Errors are computed at the final time by comparison to the exact solution. 
Figure \ref{fig:spectral} shows  the $log_{10}$ of the  $l_2$ norm of the errors vs. 
$log_{10}$ of the number of points.  The slopes of these lines (i.e the orders) are 
calculated by MATLAB's {\tt polyfit} function,
and demonstrate that the methods achieved the expected linear convergence rates.

\begin{figure}[h]
\includegraphics[scale=.5]{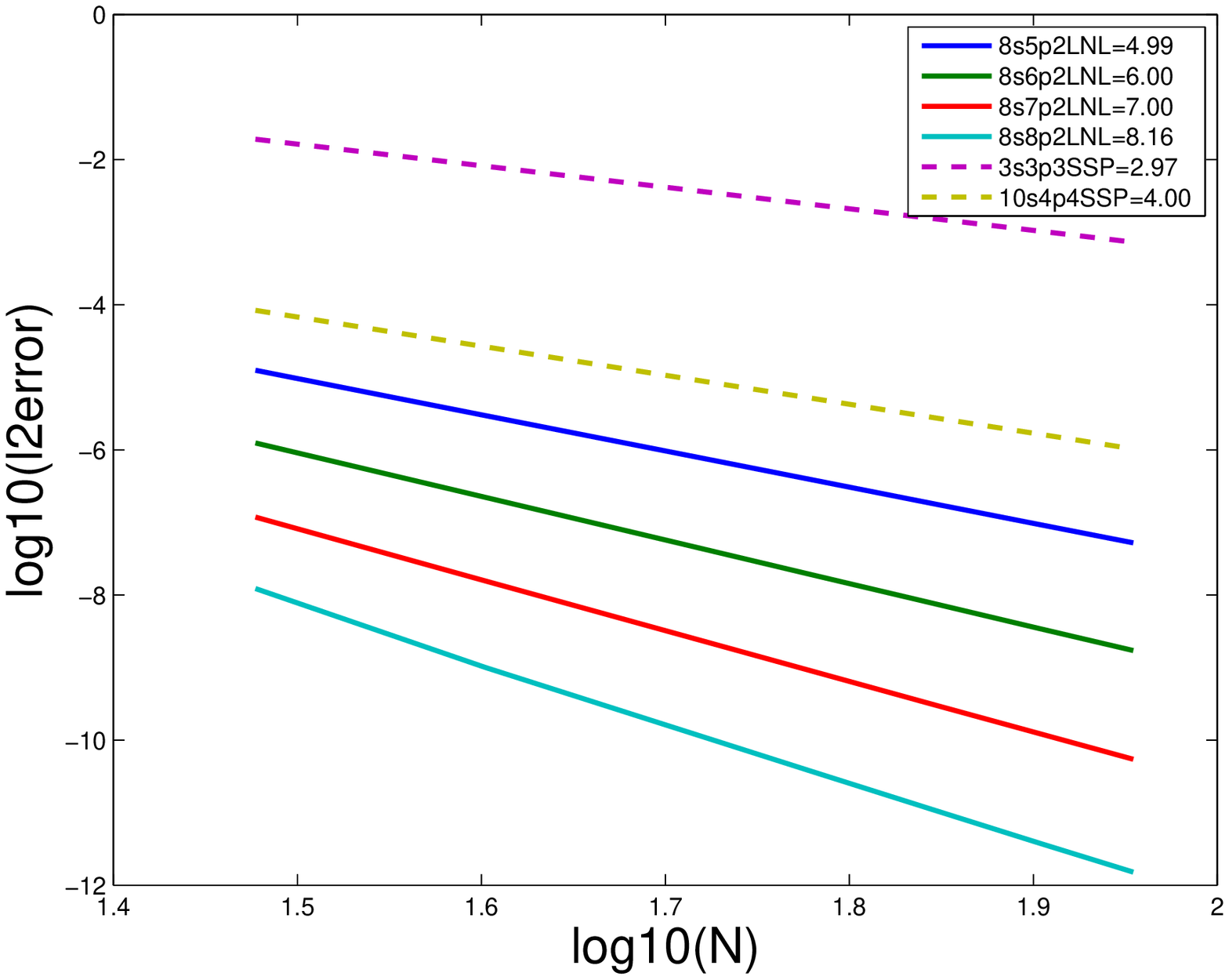}
\includegraphics[scale=.5]{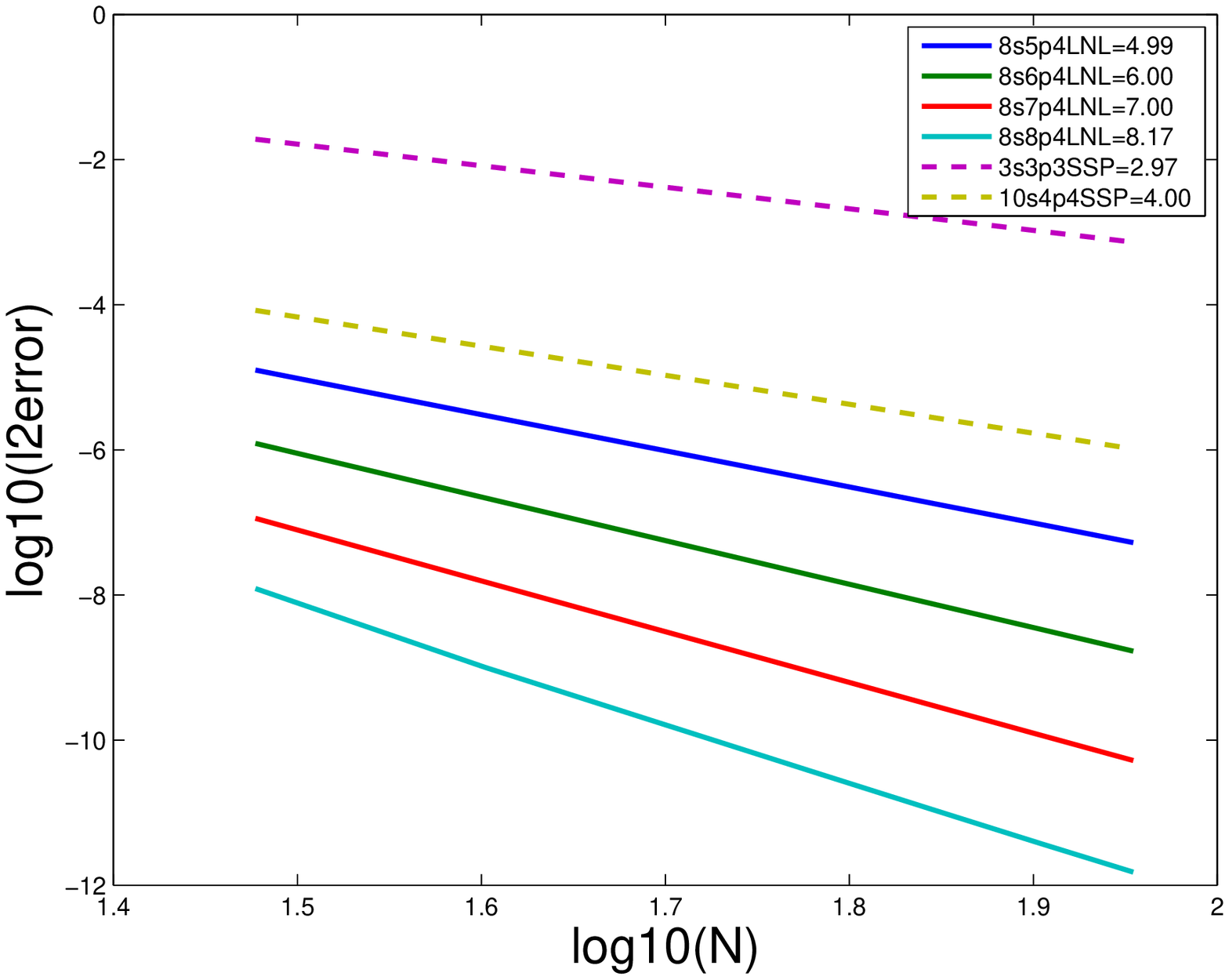}
\caption{A  linear advection problem with smooth solution and a WENO spatial discretization. 
On the left, $p=2$, on the right $p=4$. In both cases, the linear order dominates.} \label{fig:LinWENO}
\end{figure}   
\noindent{\bf Example 3: Linear advection with WENO}
We repeat the example  above, this time using the 15th order ($r=8$) WENO method to discretize in space with
$\Delta x = \frac{1}{N-1}$ for $N=30, 40, \dots , 90$. The WENO method is a nonlinear method,
so that even if the PDE is linear, the resulting system of ODEs is not. However, we can decompose the WENO
method into a linear part and a nonlinear correction term that suppresses oscillations. In theory, 
when the equation is linear and solution is smooth, the WENO method is close to linear. 
We test this problem with selected LNL Runge--Kutta time discretizations of linear order 
$p_{lin} = 5,6,7,8$ and $p=2, 3, 4$, and with the Shu-Osher SSP Runge--Kutta (3,3) 
and Ketcheson's SSP Runge--Kutta (10,4).
 As above, our final time is $T_{final}=1$, and the time step $\Delta t = 0.9 \Delta x$, where $\Delta x = \frac{1}{N-1}$. 
Errors are computed at the final time by comparison to the exact solution. 
Figure \ref{fig:LinWENO} shows  the $log_{10}$ of the  $l_2$ norm of the errors vs. 
$log_{10}$ of the number of points, and  the slopes of these lines (i.e the orders) as
calculated by MATLAB's {\tt polyfit} function.
We observe that the linear order dominates for this problem, which indicates that in regions
where the problem is primarily linear and the solution smooth, the new LNL methods with higher
linear orders could be of benefit.

\smallskip

\noindent{\bf Example 4: Burgers' equation with WENO.}
In the previous example we demonstrated the advantages of using a time-stepping method with higher $p_{lin}$
with WENO in the case of a linear, smooth problem. In this example, we 
show how a higher nonlinear order $p$ is beneficial when dealing with a nonlinear equation with possibly 
discontinuous solution.
Consider Burgers' equation with symmetric sine wave initial conditions and periodic boundaries.
\begin{eqnarray}
u_t + \left( \frac{1}{2} u^2 \right)_x  =0  \; \; \; x \in [0,1]  \\ \nonumber
u(0,x)& =& \sin(2\pi x) \; \; \; \; u(t,0)=u(t,1) . \nonumber
\end{eqnarray}
This problem develops a standing shock. We use a 
15th order WENO scheme with $N$ points in space, and 
test the LNL time-stepping methods of linear order $p_{lin} = s-1 = 5,6,7$ and nonlinear order $p=2,3,4$.
We use a time-step $\Delta t = 0.45 \Delta x$ where $ \Delta x = \frac{1}{N-1}$.

\begin{figure}[H]
\includegraphics[scale=.45]{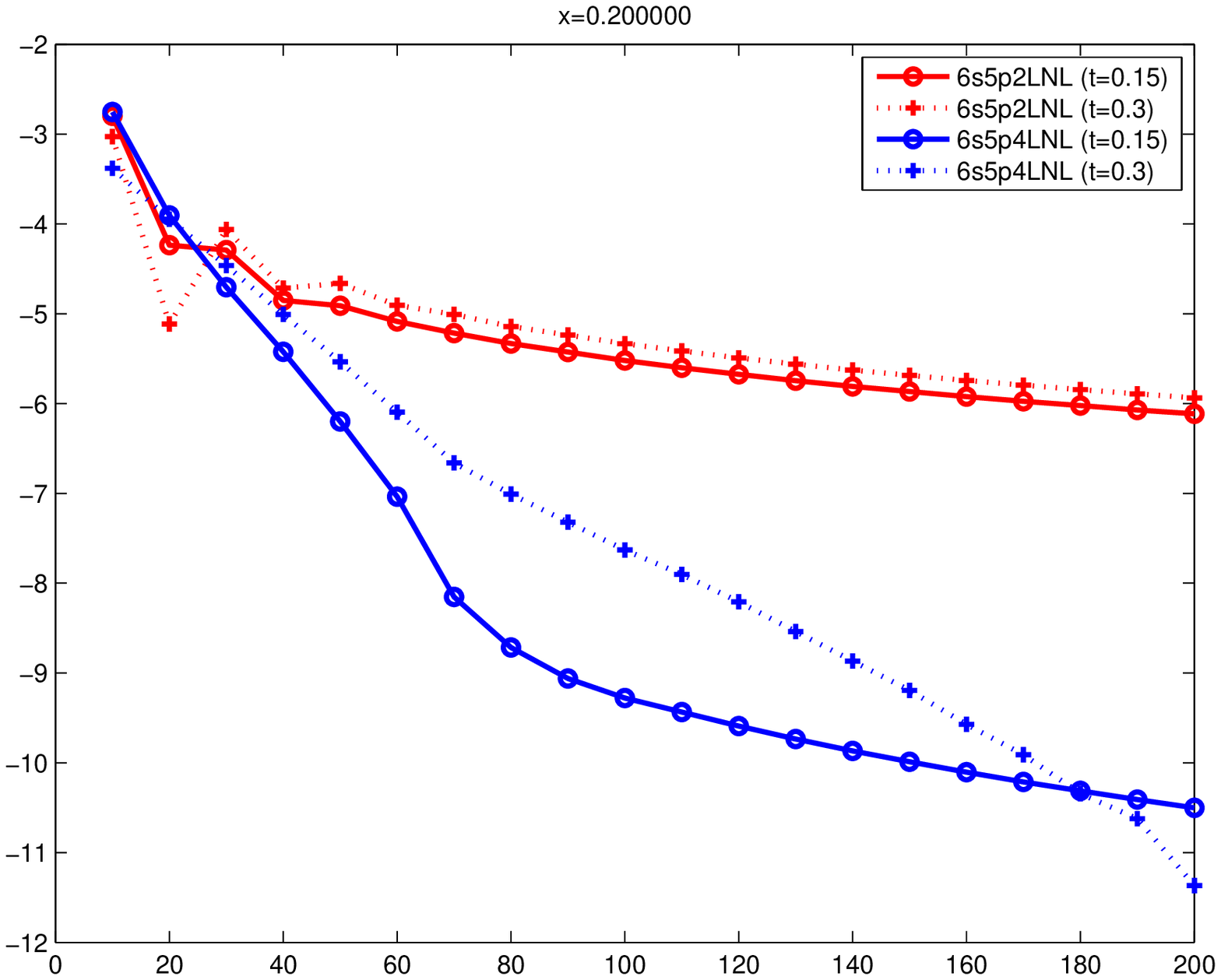} 
\includegraphics[scale=.45]{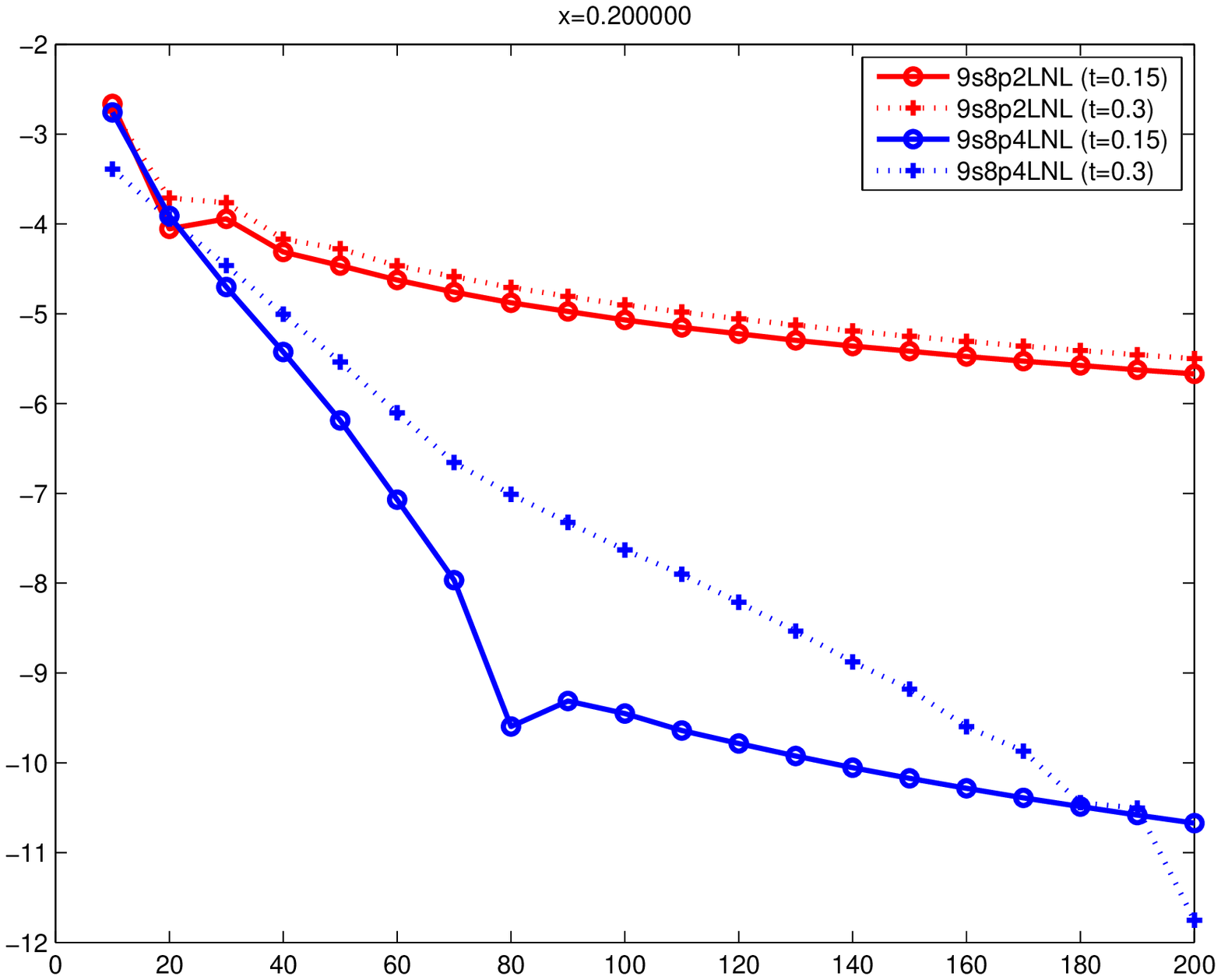}\\
\includegraphics[scale=.45]{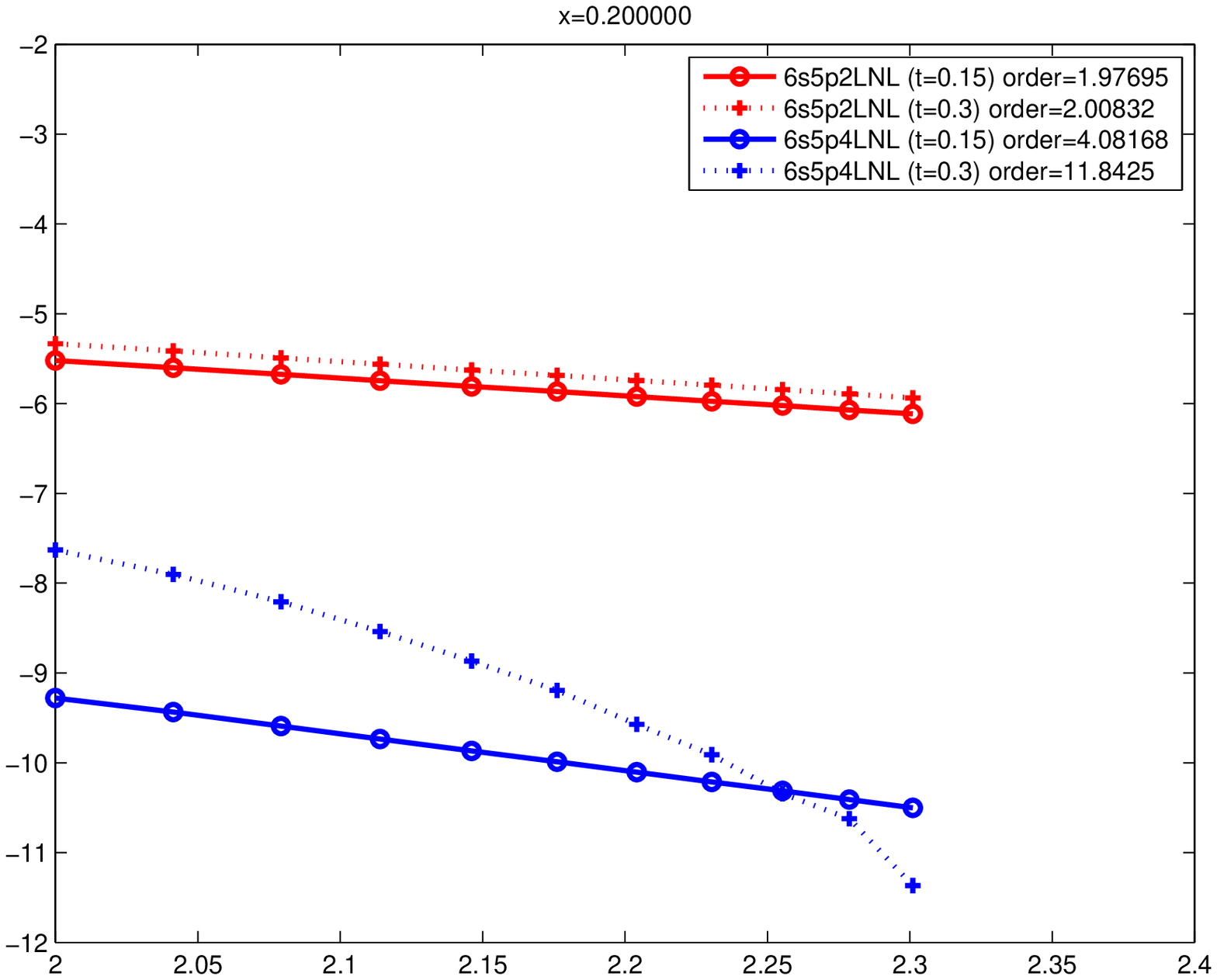} 
\includegraphics[scale=.45]{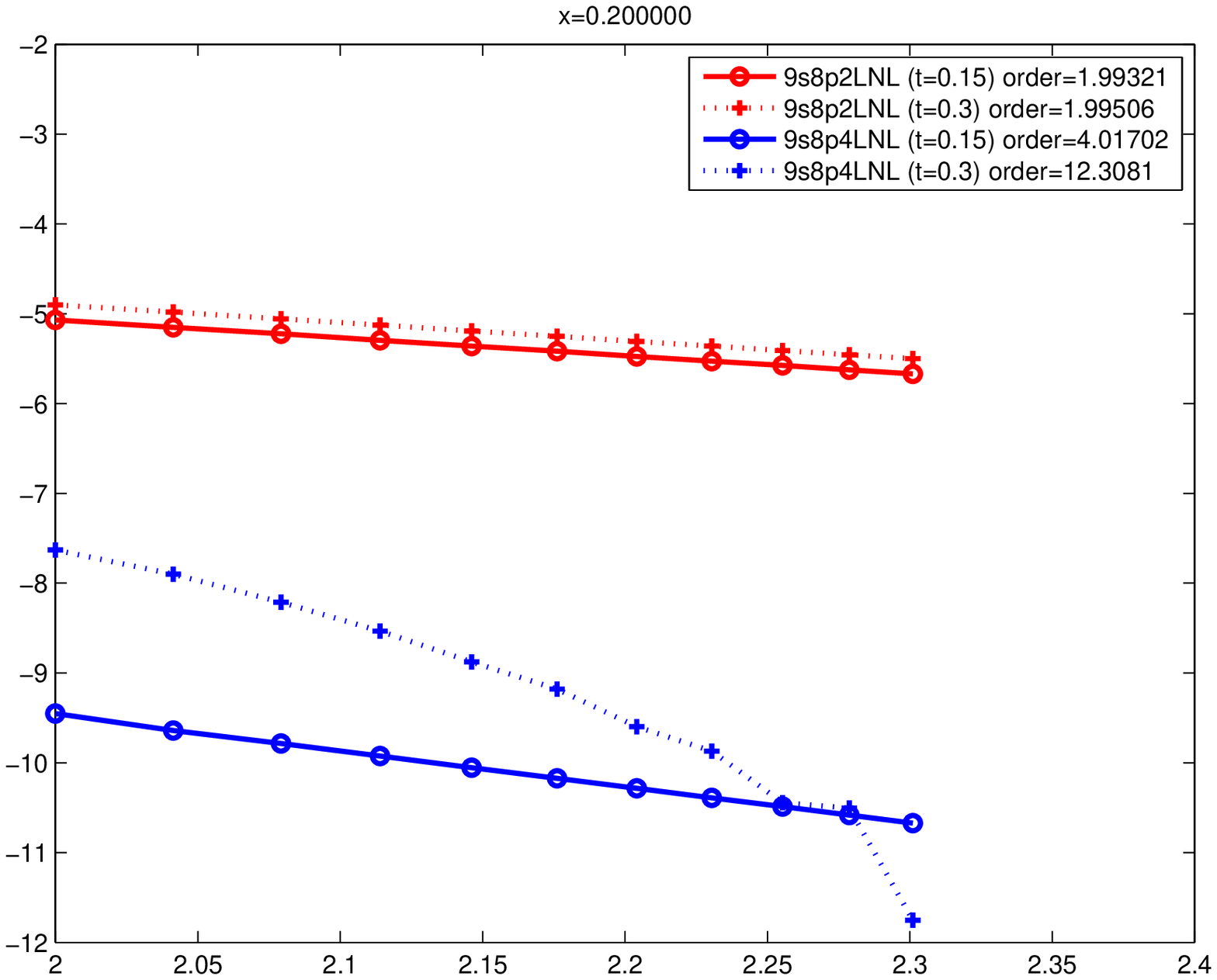}
\caption{Burgers' equation with  WENO spatial discretization. Shown are pointwise errors at $x=0.2$, before (solid lines) and 
after (dotted line) the shock. The top figures show the $log_{10}$ of the errors vs. the number of points $N$. The bottom
plots show the  $log_{10}$ of the errors vs.  $log_{10}(N)$ for $N\geq 100$.
 The time stepping methods use are $p_{lin} = s-1 = 5$ and $p_{lin} = s-1 = 8$  
with nonlinear orders  $p=2$ (red) and $p=4$ (blue). Clearly, methods with 
higher nonlinear order in time give smaller errors.}
 \label{fig:BurgersWENO}
\end{figure}

In Figure \ref{fig:BurgersWENO} we show the absolute values of the pointwise errors at spatial location 
$x=0.2$ for $N=10, \dots , 200$ (top) and for $log_{10}(N)$ (bottom). These errors are shown before the shock forms 
(at time $T_{pre} = 0.15$, solid line) and after the shock forms (at time $T_{post} = 0.3$, dotted line).
Observe that for smaller number of spatial points the errors decays very fast, however once we reach a
spatial refinement that is small enough we see that the methods with higher $p$ have significantly smaller errors.
If we consider only $N\geq 100$ points, we see the nonlinear order generally dominating: the linear methods feature second order convergence both pre- and post-shock, while the 
 $p=4$  methods are fourth order pre-shock, but jump to twelfth order post-shock (probably capturing the 
high order WENO behavior). 
Taken together with the problem in Example 3, this suggests that using a method with high linear and high nonlinear order
may be beneficial in examples that have smooth and linear regions and non-smooth nonlinear regions.

\begin{figure}[H]
 \subfigure[Fifth order methods]{\includegraphics[scale=.45]{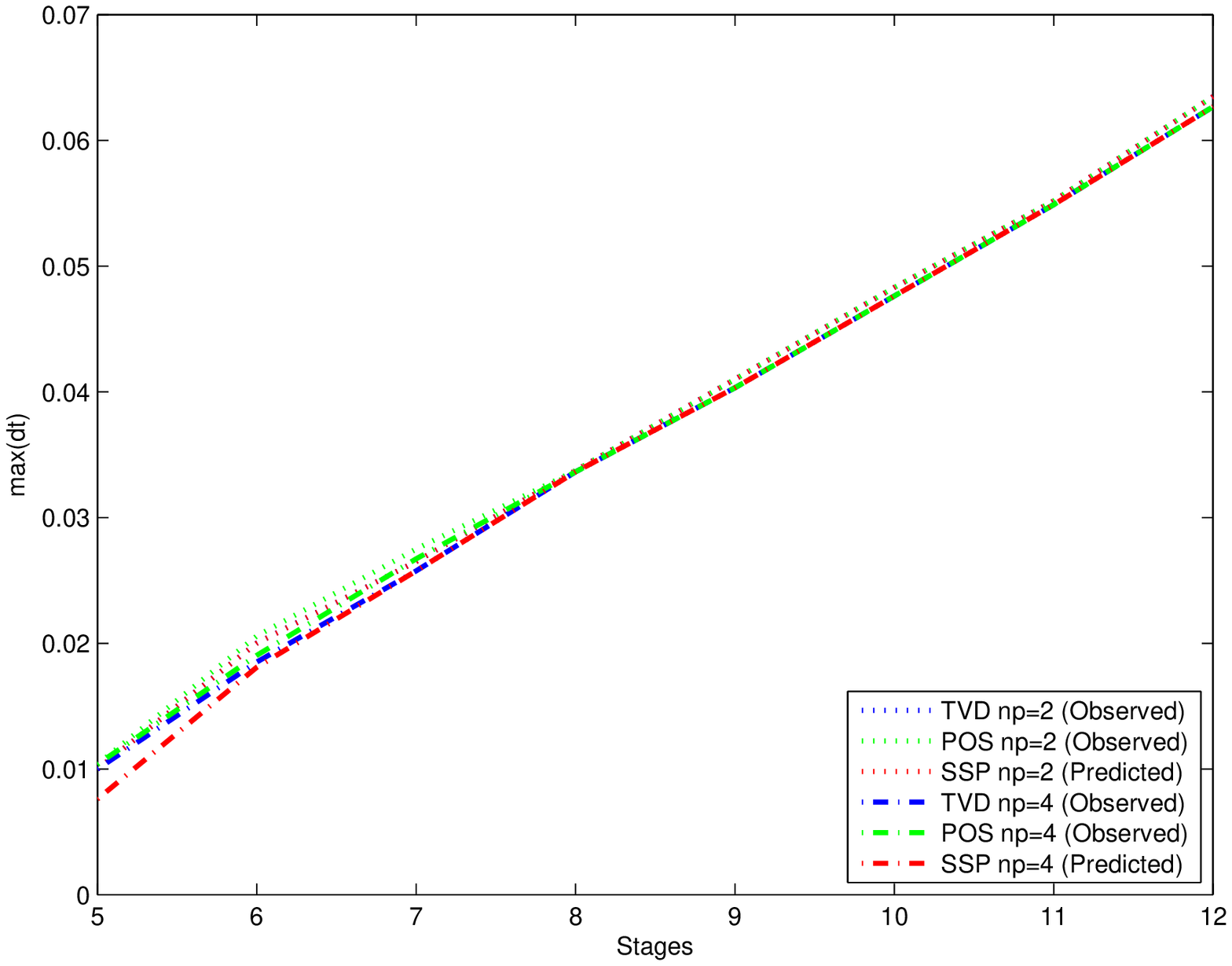}}
 \subfigure[Sixth order methods]{\includegraphics[scale=.45]{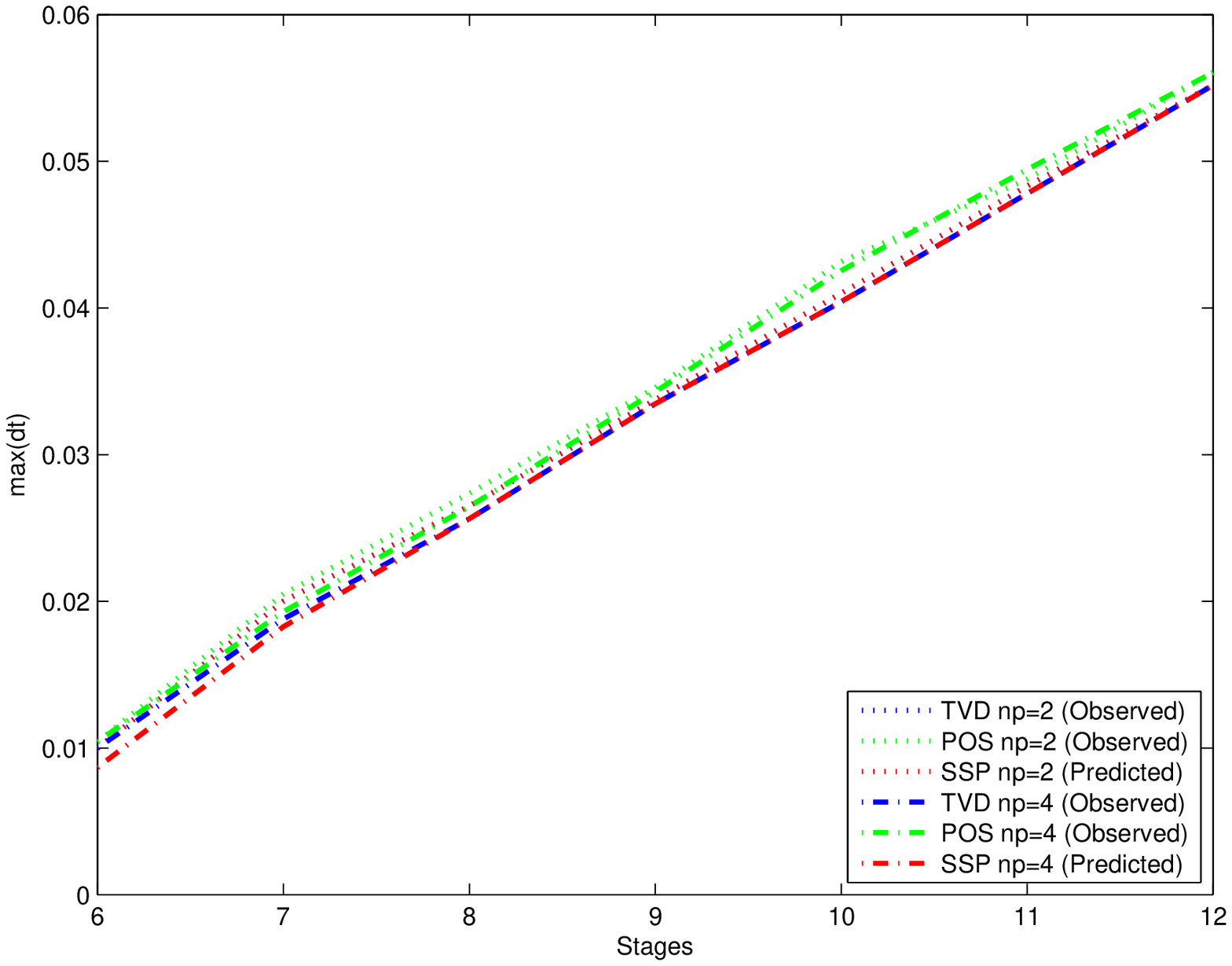}} \\
 \subfigure[Seventh order methods]{\includegraphics[scale=.45]{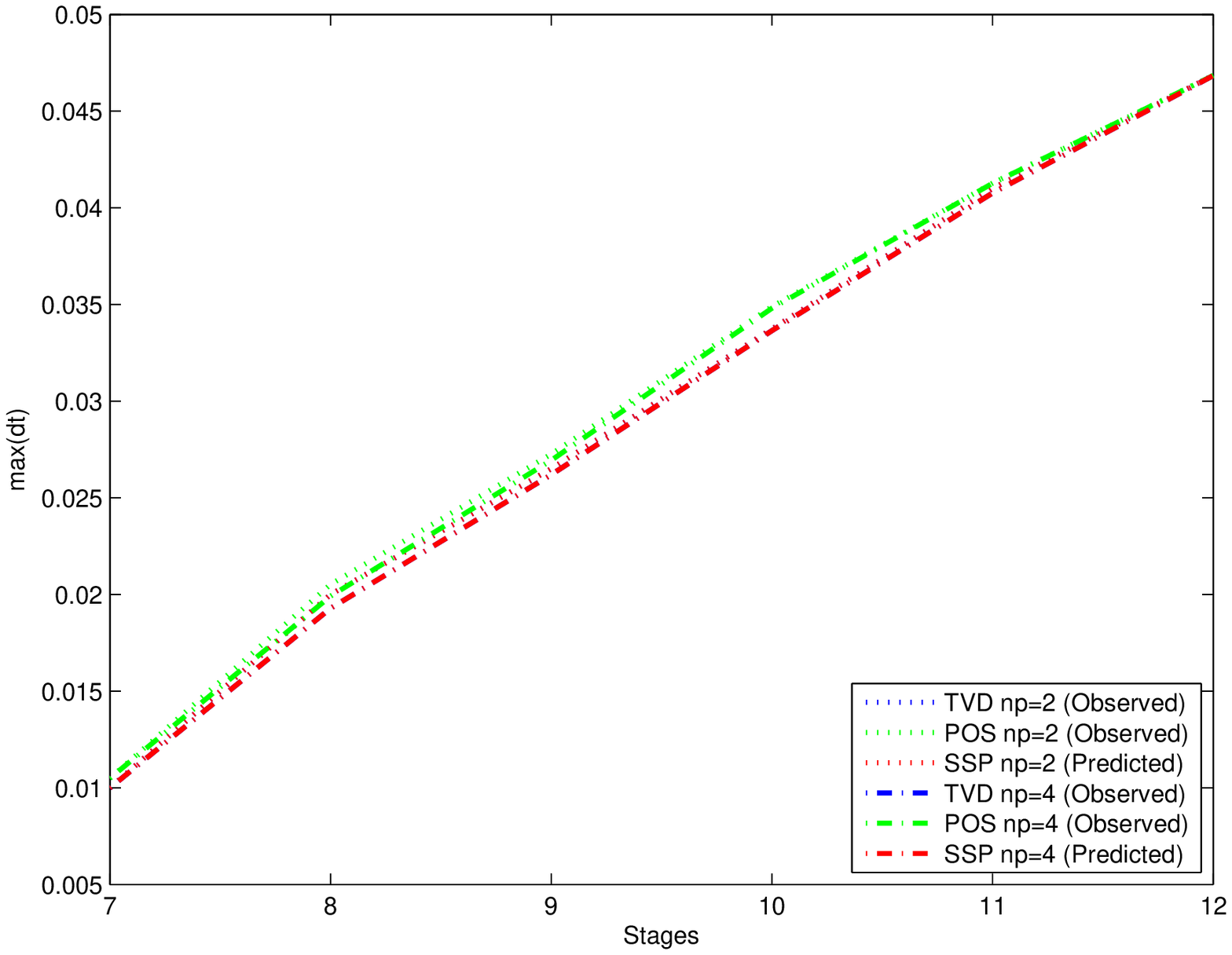}} 
 \subfigure[Eighth order methods]{\includegraphics[scale=.45]{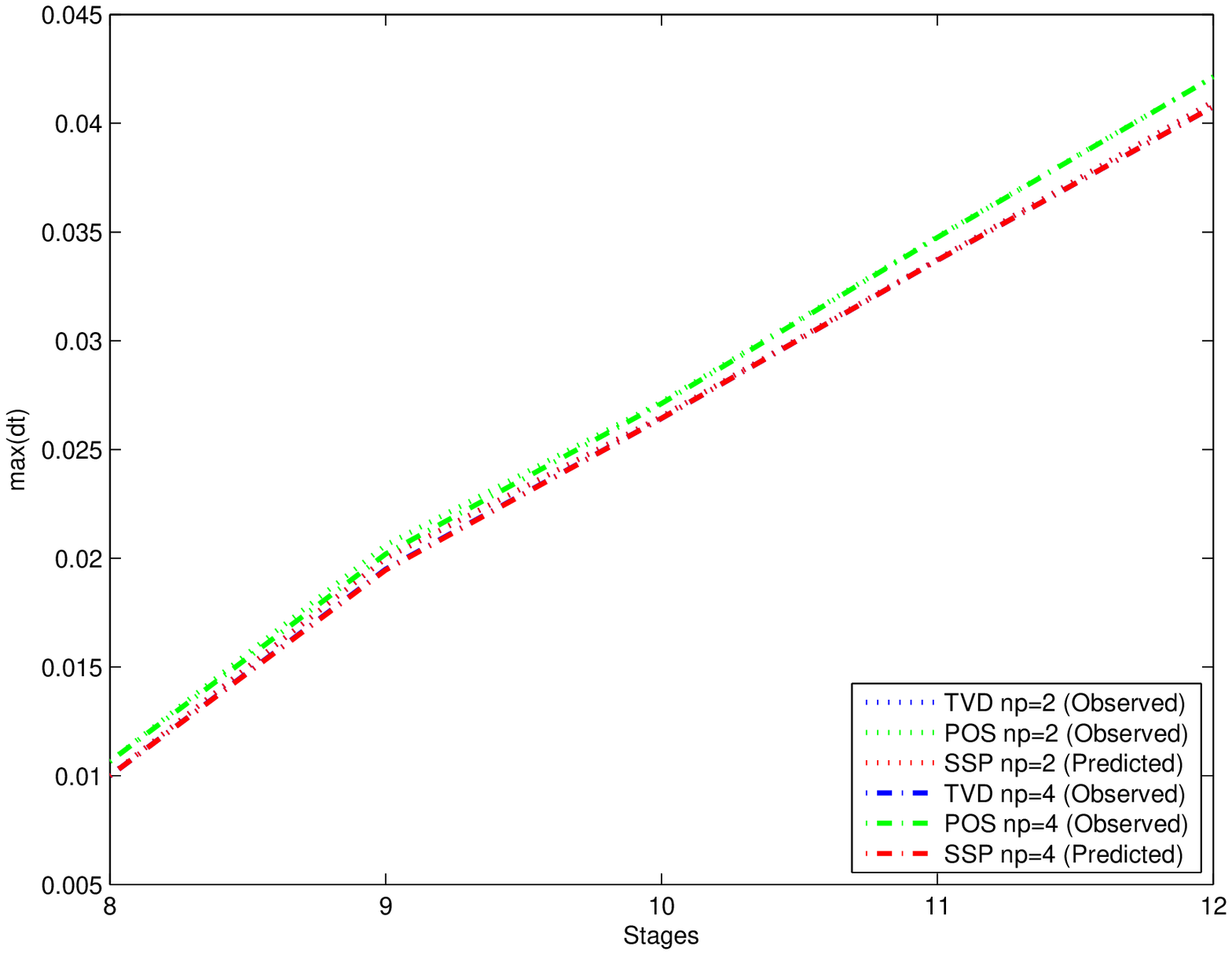}} \\
 \subfigure[Ninth order methods]{\includegraphics[scale=.45]{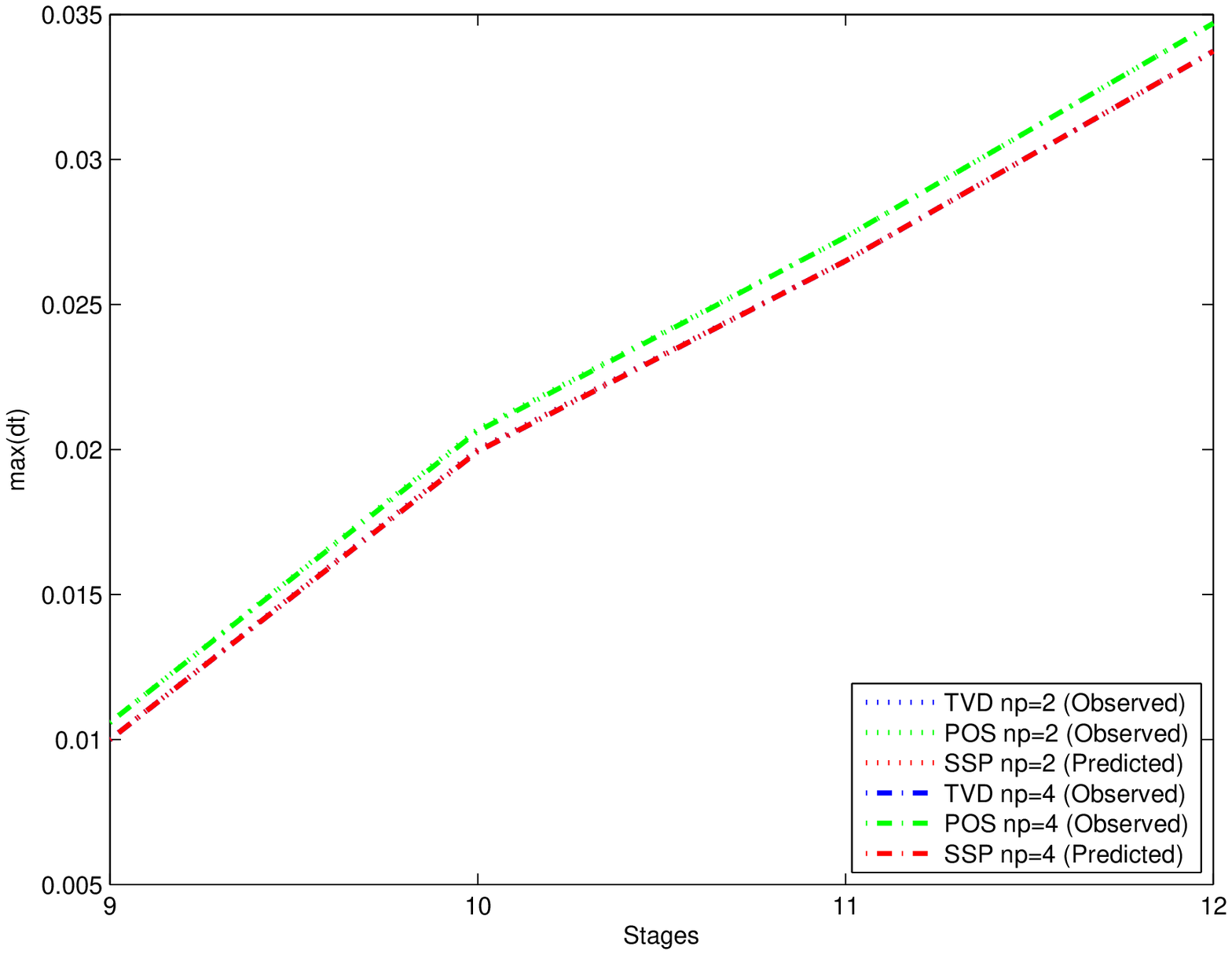}}
 \subfigure[Tenth order methods]{\includegraphics[scale=.45]{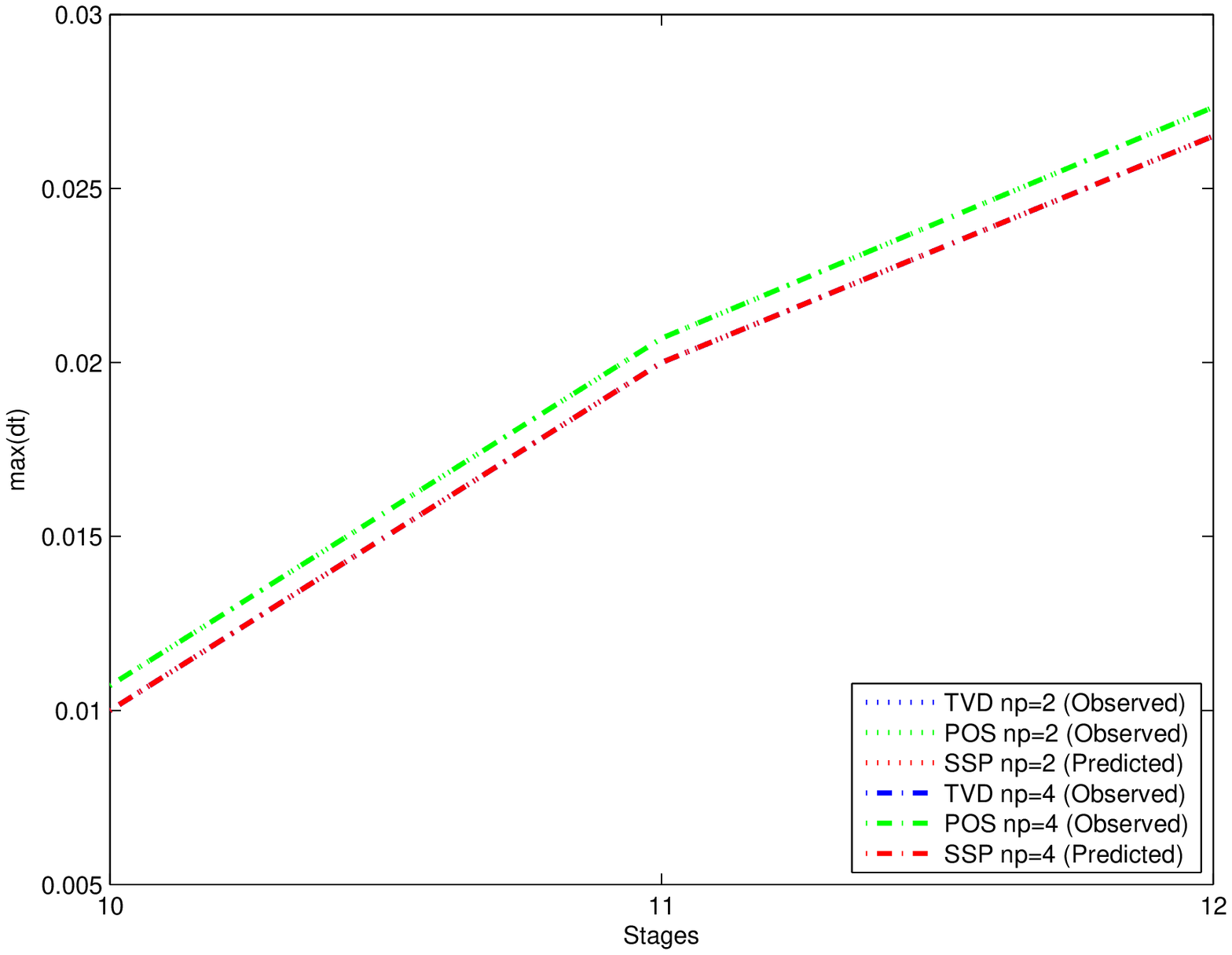}}
\caption{ The  time-step required for TVD (blue) and positivity (green) compared to the theoretical (red)
for SSP  linear (dotted) and LNL (dashed) Runge--Kutta methods of order $p_{lin}=5, . . ., 10$ 
 for a linear advection equation (Example 5).The red and blue lines overlap.}
\label{fig:LinAdvFD}
\end{figure}

\newpage
\noindent{\bf Example 5: Positivity and TVD time-step for a linear advection equation with first order finite difference
in space.} Consider the linear advection equation with a step function initial condition:
\begin{align}
u_t + u_x & = 0 \hspace{.75in}
    u(0,x)  =
\begin{cases}
1, & \text{if } 0 \leq x \leq 1/2 \\
0, & \text{if } x>1/2 \nonumber
\end{cases}
\end{align}
on the domain $[0,1]$ with periodic boundary conditions.
We used a first-order forward difference to semi-discretized this problem on a grid 
with $N=101$ points and evolved it to a final time of $t=\frac{1}{8}$. For this problem it
is known that Euler's method is TVD and positive for step sizes up to $\DtFE = \Dx$.
We computed the numerical solution using all the $p=2$ and $p=4$ SSP LNL Runge--Kutta methods
described in Section \ref{sec:optimal} and, for each one, found the largest  $\Delta t$ for which TVD and positivity
are preserved. In Figure \ref{fig:LinAdvFD} we plot these values (blue for TVD, green for positivity) compared to the 
 time-step guaranteed by the theory, $ \sspcoef \Delta t_{FE}$ (in red). 
 We observe that the theoretical value is an excellent predictor of the observed TVD time-step;
 In fact, the blue line is frequently indistinguishable from the red line. The positivity preserving time-step
 is slightly larger, but follows a similar line.

\begin{figure}[H]
 \subfigure[Fifth order methods]{\includegraphics[scale=.45]{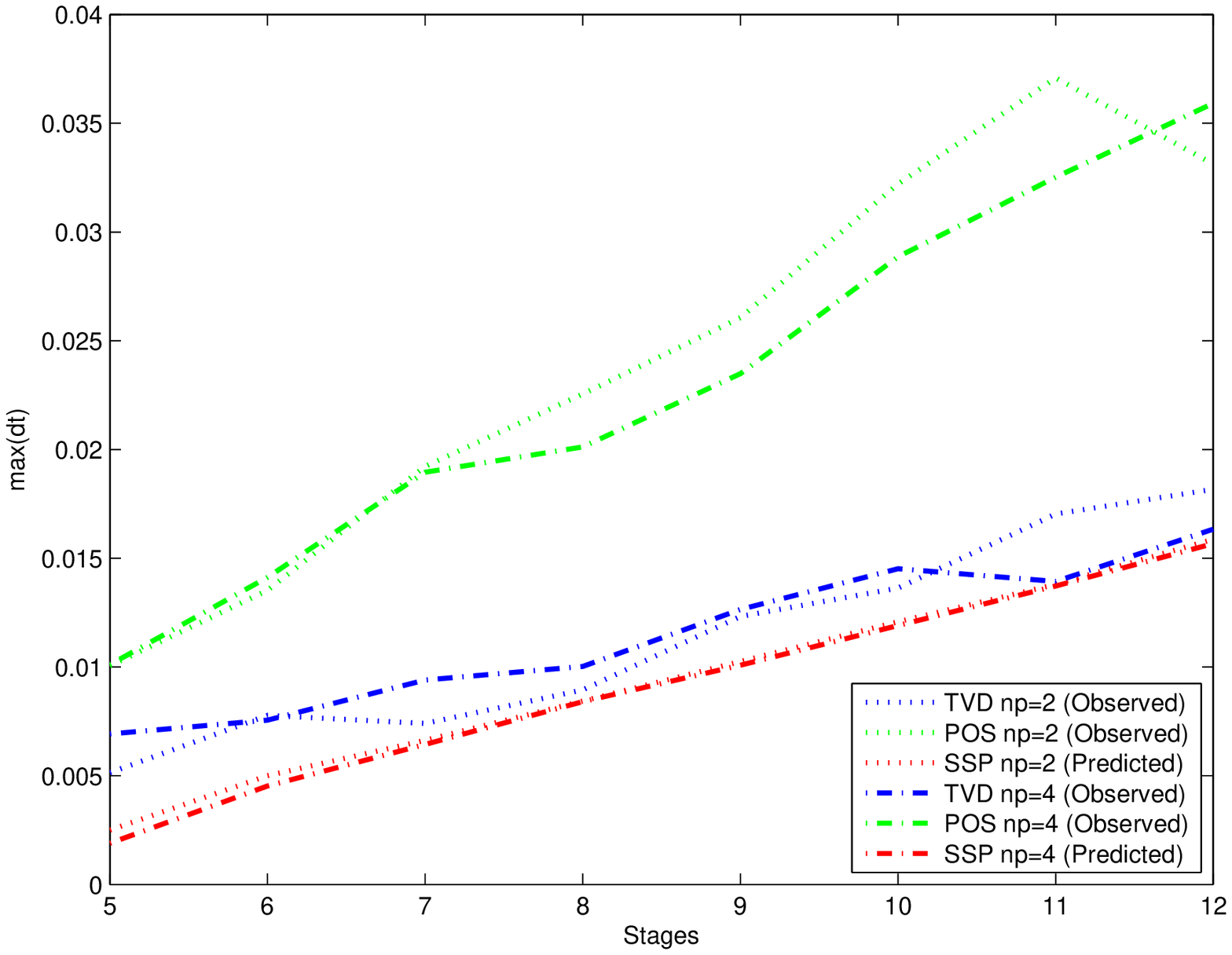}}
 \subfigure[Sixth order methods]{\includegraphics[scale=.45]{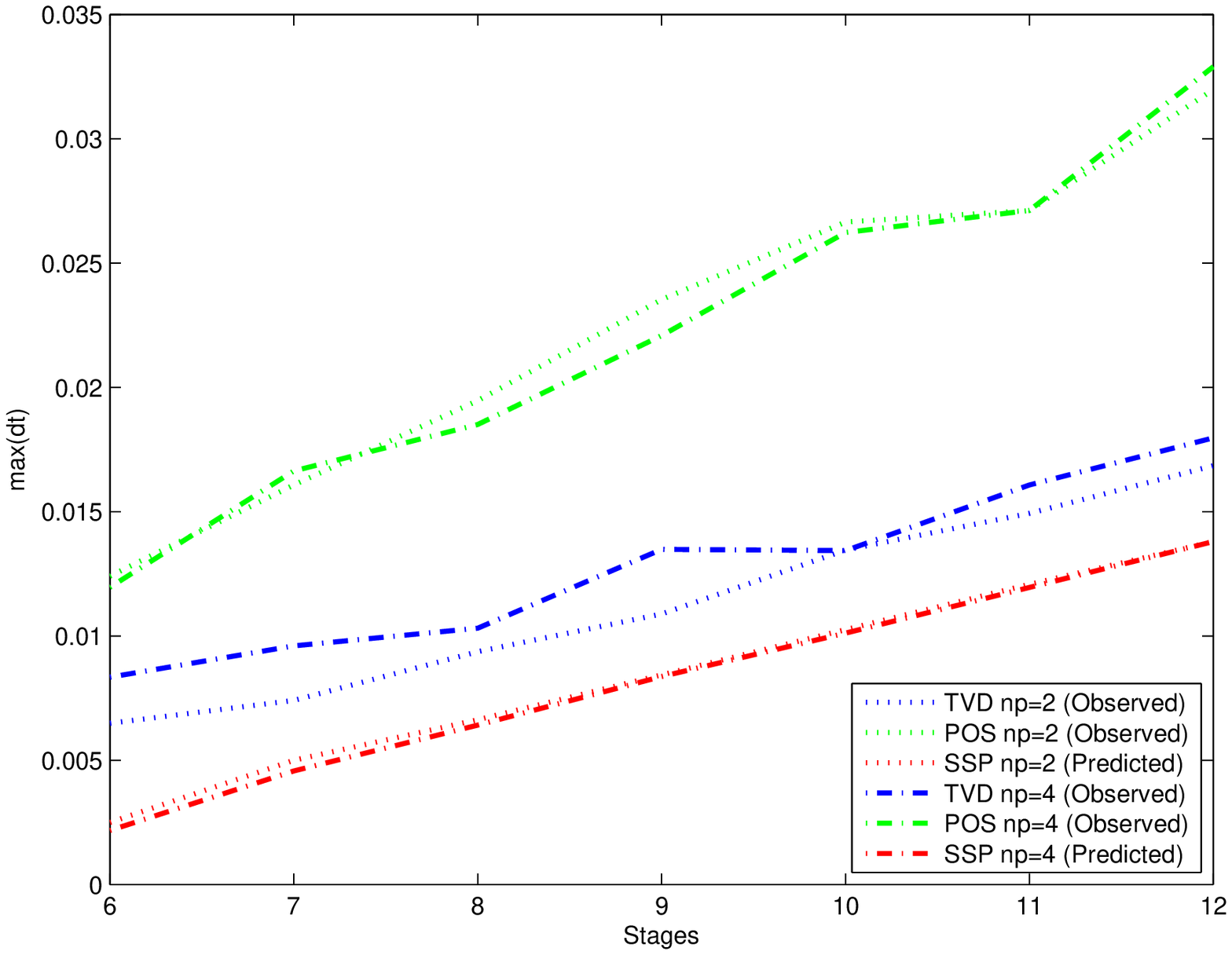}}\\
 \subfigure[Seventh order methods]{\includegraphics[scale=.45]{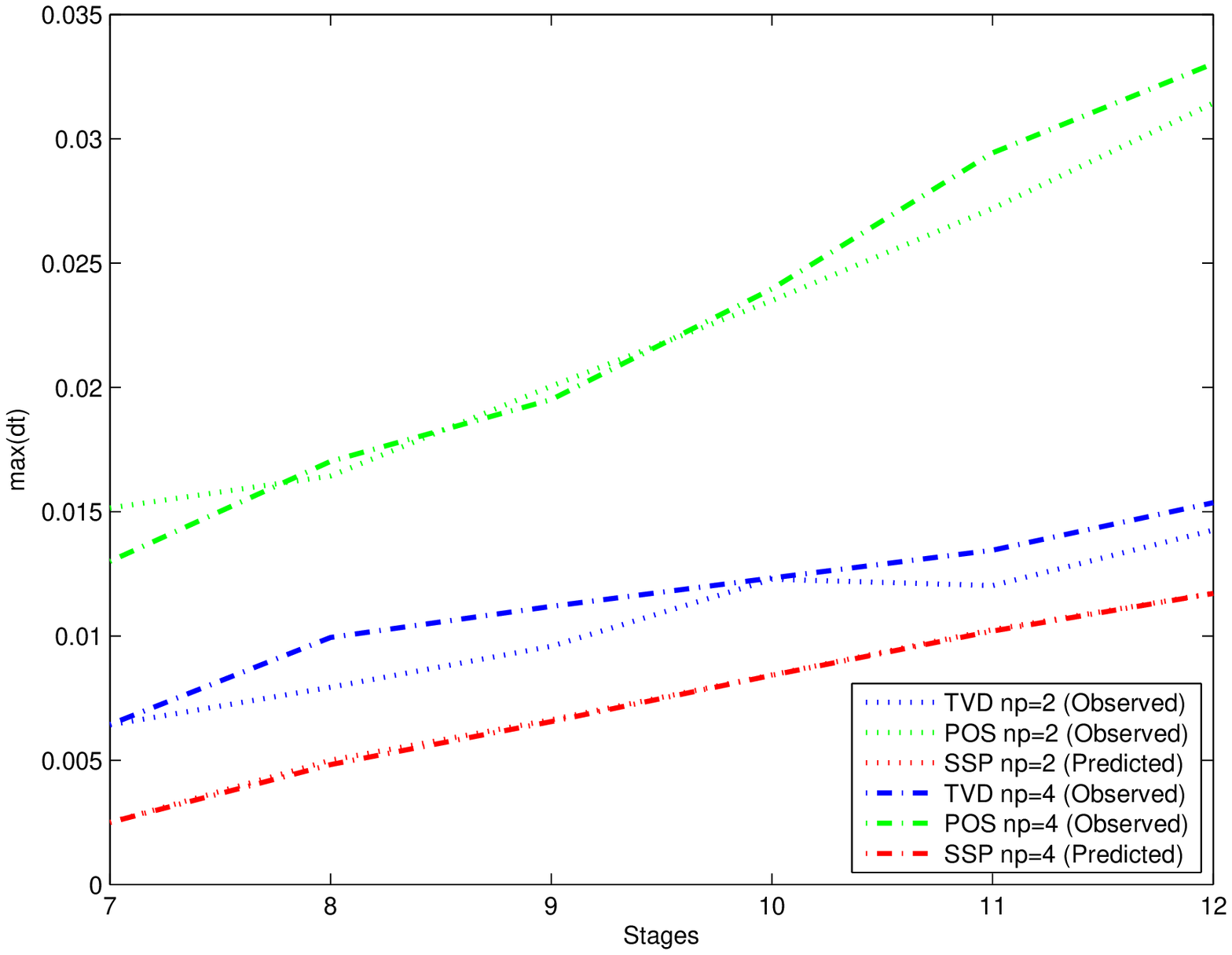}} 
 \subfigure[Eighth order methods]{\includegraphics[scale=.45]{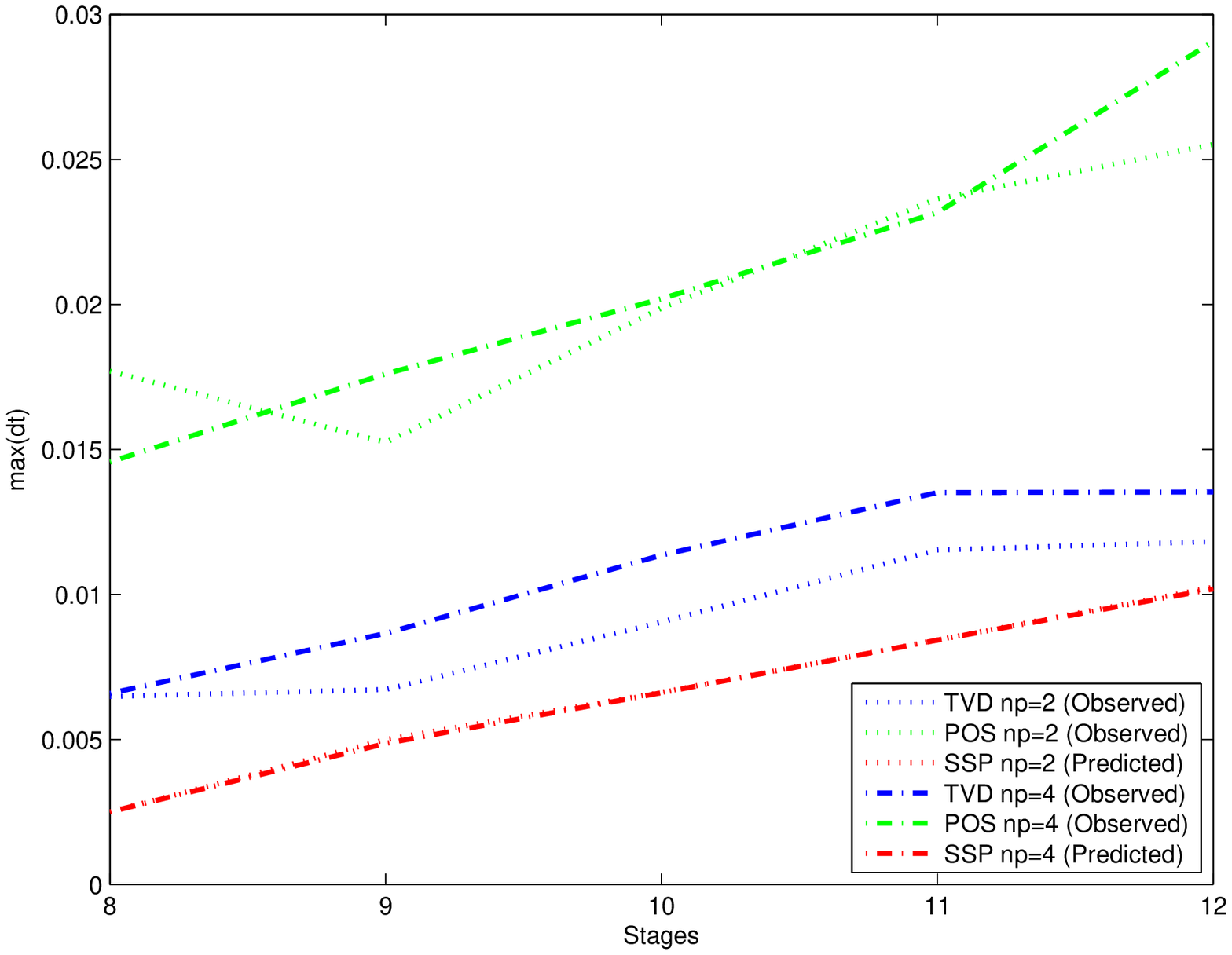}}\\
 \subfigure[Ninth order methods]{\includegraphics[scale=.45]{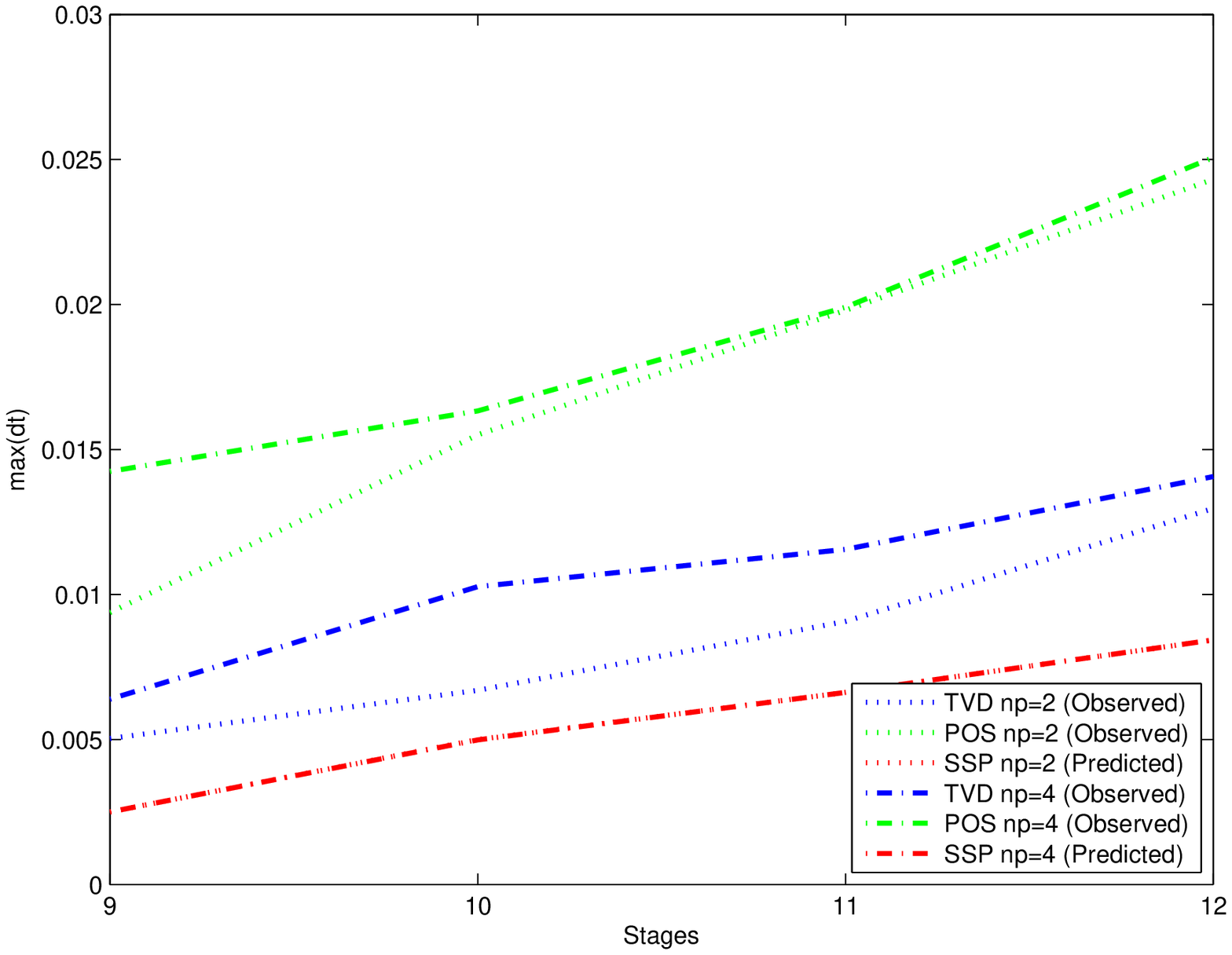}}
 \subfigure[Tenth order methods]{\includegraphics[scale=.45]{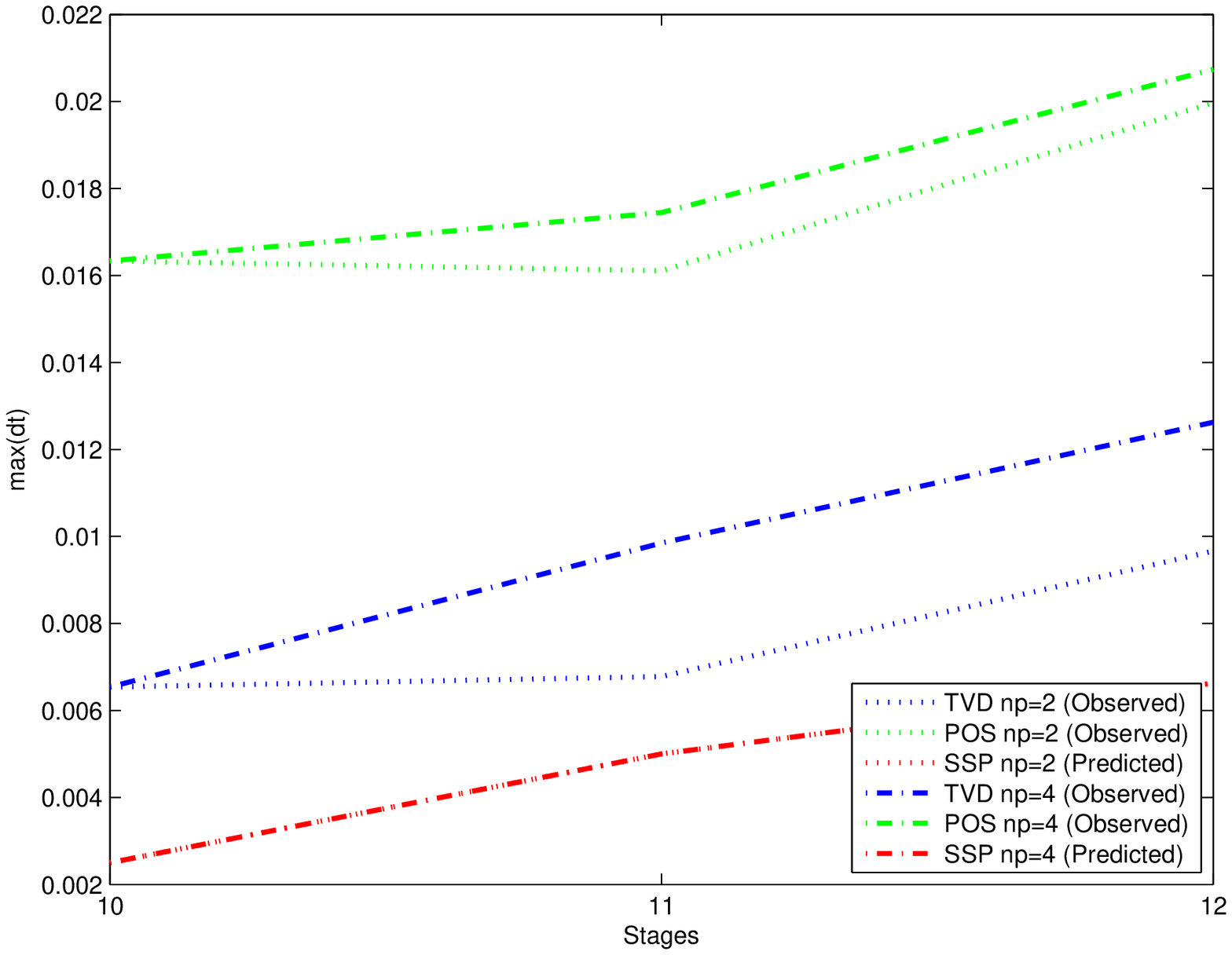}}
\caption{ The  time-step required for TVD (blue) and positivity (green) compared to the theoretical (red)
for SSP linear (dotted) and LNL (dashed) Runge--Kutta methods of order $p_{lin}=5, . . ., 10$ 
 for a nonlinear PDE (Example 6). }
\label{fig:BLtvd}
\end{figure}

\noindent{\bf Example 6: TVD and positivity for Buckley-Leverett with centered scheme and Koren limiter.}
We solve the Buckley-Leverett equation, a nonlinear PDE used to model two-phase flow through porous media:
\begin{align*}
    u_t+f(u)_x & = 0, & \text{ where } f(u) = \frac{u^2}{u^2 +a(1-u)^2},
\end{align*}
on $x\in[0,1]$, with periodic boundary conditions. 
We take $a=\frac{1}{3}$ and initial condition
\begin{eqnarray}
    u(x,0) =
\begin{cases}
1/2, & \text{if }x\ge1/2 \\
0, & \text{otherwise.}
\end{cases}
\end{eqnarray}
The problem is semi-discretized using a second order conservative scheme
with a Koren Limiter as in \cite{tsrk} with $ \Delta x = \frac{1}{100}$, and run to $t_f = \frac{1}{8}$. 
 For this problem Euler's method is TVD for $ \Delta t \leq \Delta t_{FE} = \frac{1}{4} \Delta x = 0.0025$. 
We computed the numerical solution using all the $p=2,3,4$ SSP LNL Runge--Kutta methods
described in Section \ref{sec:optimal} and,  as above found the largest  $\Delta t$ for which TVD and positivity
are preserved. In Figure \ref{fig:BLtvd} we plot these values (blue for TVD, green for positivity) compared to the 
 time-step guaranteed by the theory, $ \sspcoef \Delta t_{FE}$ (in red). 
 The observed TVD and positivity time-step are typically significantly larger than the theoretical value.
 As before, the positivity preserving time-step is larger than the TVD time-step.

\section{Conclusions}
Using the optimization procedure described in \cite{SSPbook2011,ketchcodes}, we find SSP-optimized explicit Runge--Kutta
methods that have nonlinear order of $p=3$ and $p=4$  and a higher $p_{lin}>4$ order of convergence 
on linear autonomous problems. The  order barrier of $p\leq 4$ for explicit SSP Runge--Kutta methods 
indicates the critical importance of the nonlinear order on the SSP property. Nevertheless, we find that the
{\em size} of the  SSP coefficient is typically more constrained by the linear order conditions. 
As the number of stages increases, the SSP coefficient 
becomes primarily a  function of  the relationship between the number of stages and the
linear order of the method, and not the nonlinear order. This means that in many cases, we can 
obtain methods of nonlinear order $p=4$ and linear order $p_{lin}>4$ that have the same SSP
coefficient as methods with nonlinear order $p=2$ and linear order $p_{lin}>4$. We verified the linear
and nonlinear orders of convergence of the new methods on a variety of test cases. We also showed
the behavior of these new LNL time-stepping methods coupled with the WENO method for both linear and
nonlinear problems, which suggests that these LNL methods may be useful for problems that have regions
that are dominated by linear, smooth solutions and other regions where the solution is discontinuous
or dominated by nonlinear behavior. Finally, we studied the total variation diminishing and positivity preserving
properties of these LNL methods on linear and nonlinear problems, and showed that for the linear problems, the
theoretical SSP time-step is a very accurate predictor of the observed behavior, while serving only as a lower
bound in the nonlinear case. We conclude that where methods with high linear order are desirable,
it is usually advantageous to pick those methods that also have higher nonlinear order ($p>2$).

{\bf Acknowledgment.} The authors wish to thank Prof. Bram van Leer
for the motivation for studying this problem, and Prof. David Ketcheson for many helpful discussions.
This publication is based on work supported by  AFOSR grant FA-9550-12-1-0224
and KAUST grant FIC/2010/05.

\end{document}